\documentclass{article}
\usepackage{lipsum}
\usepackage{amsfonts}
\usepackage{graphicx}
\usepackage{epstopdf}
\usepackage{amsmath,amssymb,bm, cancel,xcolor}
\usepackage{amsthm}
\usepackage{algorithmic}
\usepackage{hyphenat}
\usepackage{mathrsfs}
\newcommand{\norm}[1]{\left\lVert#1\right\rVert}
\usepackage{amsthm}

\makeatletter
\def\thmheadbrackets#1#2#3{
\thmname{#1}\thmnumber{\@ifnotempty{#1}{ }\@upn{#2}}
\thmnote{ {\the\thm@notefont[#3]}}}
\makeatother

\newtheoremstyle{brakets}
{}
{}
{\itshape}
{}
{\bfseries}
{.}
{ }
{\thmheadbrackets{#1}{#2}{#3}}

\newtheoremstyle{defbrakets}
{}
{}
{\normalfont}
{}
{\bfseries}
{.}
{ }
{\thmheadbrackets{#1}{#2}{#3}}

\theoremstyle{brakets}
\theoremstyle{defbrakets}

\newtheorem{theorem}{\bf Theorem}[section]
\newtheorem{prop}{\bf Proposition}[section]

\newtheorem{cor}{\bf Corollary}[section]
\newtheorem{definition}{\bf Definition}[section]
\newtheorem{lemma}{\bf Lemma}[section]
\newtheorem{remark}{\bf Remark}[section]

\newcommand{\uu}{{\mathbf u}}

\renewcommand{\a}{\mathbf{n}}

\newcommand{\p}{\mathbf{p}}
\newcommand{\q}{\mathbf{q}}

\renewcommand{\S}{Sec. }

\newcommand{\ee}{\mathbf{e}}

\usepackage[colorlinks=true,linkcolor=blue,citecolor=blue]{hyperref}
\usepackage[
backend=biber,
style=alphabetic,
sorting=ynt
]{biblatex}
\renewenvironment{proof}{{\bfseries Proof}}{\qed}
\bibliography{mybib}
\usepackage{authblk}

\begin{document}

\title{Generalizing {Parametrization} Invariance in the Calculus of Variations}

\author[a]{Sanjay Dharmavaram\thanks{sd045@bucknell.edu}}
\author[b]{Basant Lal Sharma\thanks{bls@iitk.ac.in}}

\affil[a]{Department of Mathematics, Bucknell University, 17837 PA, U.S.A.}
\affil[b]{Department of Mechanical Engineering, Indian Institute of Technology Kanpur, Kanpur, 208016 UP, India}

\maketitle

\begin{abstract}
We revisit the notion of {parametrization} invariance while introducing certain weakened notions of invariance in the calculus of variations. In this work, we {employ a straightforward approach in} the classical setting and mostly restrict attention to functionals on one-dimensional domains. We establish a connection between {parametrization} invariant functionals and functionals embodying a weaker notion of invariance of their Lagrangian; we term this notion as \emph{{$\mathcal{T}$-}{}Lagrangian} analogous to the well-known idea of null Lagrangian. However, the Euler-Lagrange operator of a {$\mathcal{T}$-}{}Lagrangian vanishes only along the tangential direction {in the configuration space}. On one-dimensional domain {and for first- and second-order theories}, we show that functional described by such a Lagrangian is {necessarily} a {parametrization} invariant functional modulo null Lagrangian. Keeping the motivation for partial differential equations, we also introduce and explore the notion of \emph{{$\mathcal{N}$-}{}Lagrangian}, {with an invariance complementary to the case of \emph{{$\mathcal{T}$-}{}Lagrangian}}, whose Euler-Lagrange operator {vanishes} along normal directions. 
We find that in a one-dimensional setting, every {$\mathcal{N}$-}{}Lagrangian is simply a null Lagrangian.
\end{abstract}


\section*{Introduction}
{Parametrization} invariance, {often called parametric invariance and} also sometimes called reparametrization invariance, is related to an important symmetry of many physical theories \cite{gueorguiev2021reparametrization} {due to their underlying geometric structure}. {For instance,} it manifests in the theory of general relativity \cite{thorne2000gravitation} as a consequence of the principle of general covariance \cite{norton1993general} -- that the equations of physics must be independent of the choice of coordinate systems. Most theories of quantum gravity, including string theory, therefore, inherit this symmetry \cite{horowitz1986origin}.

In physics, {parametrization} invariance is also encountered as an important example of gauge symmetry \cite{fulop1999reparametrization}, a symmetry characterized by an infinite-dimensional Lie group of transformations of a functional ({often referred as the} action). {The} symmetry group for {parametrization} invariance is {a} diffeomorphism group. A defining feature of \emph{gauge theories} (theories possessing gauge symmetry) is that their Euler-Lagrange equations are necessarily under-determined. This can be seen to arise as a consequence of Noether's second theorem \cite{olver1986noether, hydon2011extensions}.
In the case of {parametrization} invariant functionals, { in fact}, the component of the Euler-Lagrange operator in the tangential direction to the field variable is trivially zero.

{Parametrization} invariance is not limited to theories of high-energy physics; it also occurs with a firm footing in traditional continuum mechanics. For example, it is a defining symmetry of the Helfrich-Canham model (discussed in more detail in \S\ref{sec:motivation}) \cite{canham1970minimum, helfrich1973elastic} for lipid bilayer membranes and of the theory of minimal surfaces (a model for soap films) \cite{GiaquintaV1} {besides being an important geometric feature of the Euler elastica \cite{euler1952methodus}}. It is evident that {these} theories are independent of the parametrization chosen to describe the surface. In the Helfrich-Canham model, {parametrization} invariance accounts for the in-plane fluidity of the lipid membrane \cite{capovilla2002stresses}; the lipid molecules freely flow on the surface of the membrane, encountering very little resistance from the surrounding molecules. In this way, the {parametrization} invariance of the Helfrich-Canham energy can be ascribed to the lack of a well-defined in-plane reference configuration for the membrane. 
{On the other hand, the parametrization invariance of the problem of Euler elastica arises solely out of the elastic property of the deformation where the stored energy per unit length in deformed configuration depends only on the current curvature}.
The equations of motion (Euler-Lagrange equations) for these theories are typically derived by taking variations of the surface or to the curve in the normal and tangential directions to the surface \cite{capovilla2002stresses, colding2011course}. It can be shown \cite{steigmann2003variational} that the tangential component of the Euler-Lagrange operator is trivially zero. 

Indeed, the equations of motion{, or evolution equations, for several phenomena in physics and mechanics, including those alluded to} in previous paragraph, are typically derived {as Euler-Lagrange equations} by taking variations of the {field} in the normal and tangential directions {in configuration space considered to be embedded in a containing Euclidean space,
wherein the} starting point for {the mathematical formulation is the proposal of} an energy functional that satisfies certain symmetries ({for example, the parametrization} invariance for the examples noted above). In continuum mechanics, {particularly in the theory of elastic media}, such models are called \emph{hyperelastic models} and are {examples of} conservative systems. However, in some cases, the theory {for a physical system} is formulated by starting with the equations of motion using balance laws {such as balance of linear and angular momentum, mass, etc}; {one such} example is Nagdhi's formulation for shell theory \cite{naghdi1973theory}. In such a class of physical models, where the analog of {parametrization} invariance is provided at the level of equations of motion, it is unclear which class of functionals fulfills that role {in an exhaustive manner}. {It} is also important to note that the way to analyze models with a restricted kind of fluid-invariance is speculative too; the latter is anticipated to be relevant with the inclusion of microscopic effects or environmental changes, c.f., {chapter 10, \cite{nelson2004statistical}},\cite{Deserno}.

At this point a question that naturally arises {is} --- \emph{What is the analog of {parametrization} invariance of integrals/functionals at the level of differential equations of motion?} In this paper, we explore this question by introducing a weakened notion of invariance of a Lagrangian which we term \emph{{$\mathcal{T}$-}{}Lagrangian}. We address this issue in one dimension where we use \emph{elementary calculus} to show that a weakened notion of invariance is interesting. We prove general representation result (Theorems \ref{thmatangt_full ndof} and \ref{thm:rep thm second order}) for the {$\mathcal{T}$-}{}Lagrangian and connect it to {parametrization} invariance. We also explore the complementary notion of an {$\mathcal{N}$-}{}Lagrangian. Surprisingly, we find that (Theorem \ref{thm:Nnull=nullLag}) in the one-dimensional setting, {$\mathcal{N}$-}{}Lagrangians are exactly classical null Lagrangians \cite{dE62,Ericksen62,Rund66}. We conjecture that generalizing these concepts to higher dimensions may offer nontrivial implications. It is also plausible that {$\mathcal{N}$-}{}Lagrangian may not be as trivial as it is in the one-dimensional case; {it is worth a note that the characterization of null Lagrangians itself is a daunting task in higher dimensions and complex field theories, see \cite{sharma2021null} for a recent application of the results of \cite{OS88}. Altogether these weakened notions of invariance are anticipated to involve a significant role played by null Lagrangians. Regarding the ways of the present paper, the simplicity of presented one-dimensional setting also allows a possible analysis on lines of \cite{bates2016elastica} that exploits jet bundle formalism; it is expected that the same formalism becomes handy and has strong potential in higher dimensions where our approach becomes tedious}.

Overall, the intended plan {in this paper} concerns addressing the rephrased question: What are the functionals that satisfy the weakened notion of invariance with their Euler-Lagrange equations becoming trivial in certain directions (depending on the field), or, equivalently, remaining non-vanishing in certain directions (for some fields)? In this paper, we take only a first step towards posing this question and answering it in a simple setting while the future plan is attaining the ambitious goal of extending the introduced notions and concomitant analysis to partial differential equations.

This paper is organized as follows. In \S\ref{sec:motivation}, we motivate our work using {few} examples of {parametrization} invariant functionals. 
In \S\ref{sec_revisit}, we revisit the notion of {parametrization} invariance and briefly dwell upon a well-known property of {parametrization} invariant functionals. In \S\ref{sec_tangentvar}, we introduce the notion of {$\mathcal{T}$-}{}Lagrangian and prove the representation theorem (Theorem \ref{thmatangt_full ndof}) for such Lagrangians of {first}-order. The counterpart for second-order theories (Theorem \ref{thm:rep thm second order}) is discussed in \S\ref{sec:second order}. In \S\ref{sec_normalvar}, we introduce the notion of {$\mathcal{N}$-}{}Lagrangian and similarly derive a representation result {but only for the case of first order}. {In \S\ref{sec_discuss} we discuss few more aspects of the analysis in the context of some continuum theories.} We summarize our conclusions and point out a few related 
questions in \S\ref{sec:conclusions}.

\section{Motivation}
\label{sec:motivation}
{First we discuss} {some} examples from continuum mechanics to motivate the notion of {$\mathcal{T}$-}{}Lagrangian. These examples serve as templates for first-order and second-order theories, respectively. As already referred in the introduction to this paper, we consider: theory of minimal surfaces, the Helfrich-Canham theory for fluid membranes, {and the Euler elastica}.

A minimal surface $\omega\subset \mathbb{R}^3$ is a minimizer of its total (surface) area functional \cite{colding2011course}: 
\begin{equation}
\mathcal{E}_{MS} = \int_\omega\;da,
\label{mintheory}
\end{equation}
subject to some prescribed boundary conditions. Here $da$ is the area measure on $\omega$. The classical Plateau's problem concerns establishing the existence of such surfaces, which plays an important role in understanding the shape of soap films \cite{almgren1976geometry}.

The second example we consider is the Helfrich-Canham theory \cite{canham1970minimum, helfrich1973elastic} for modeling the bending elasticity of lipid bilayer membranes. The energy functional for the corresponding model is given by
\begin{equation}
\mathcal{E}_{HC} = \int_\omega (\kappa H^2 + \kappa_g K)\;da,
\label{HCenergy}
\end{equation}
where $\omega$ is the surface describing the deformed membrane, $H$ is the mean curvature of the surface, $K$ is the Gaussian curvature, and $\kappa$ and $\kappa_g$ are the bending moduli. {T}hese functionals depend only on geometric entities associated with $\omega$, it is independent of the parametrization chosen to describe the surface. Indeed,
{parametrization} invariance manifests to account for the in-plane fluidity of the lipid membranes \cite{capovilla2002stresses}; the lipid molecules freely flow on the surface of the membrane, encountering very little resistance from the surrounding molecules. In this way, {parametrization} invariance of the Helfrich-Canham energy \eqref{HCenergy} can be ascribed to the lack of a well-defined in-plane reference configuration for the membrane.\footnote{In the mathematical literature dealing with differential geometry, \eqref{HCenergy} is often addressed as the Willmore functional \cite{kuwert2012willmore}.}

\begin{figure}
\centering
\begin{tabular}{cc}
\includegraphics[width=.35\textwidth]{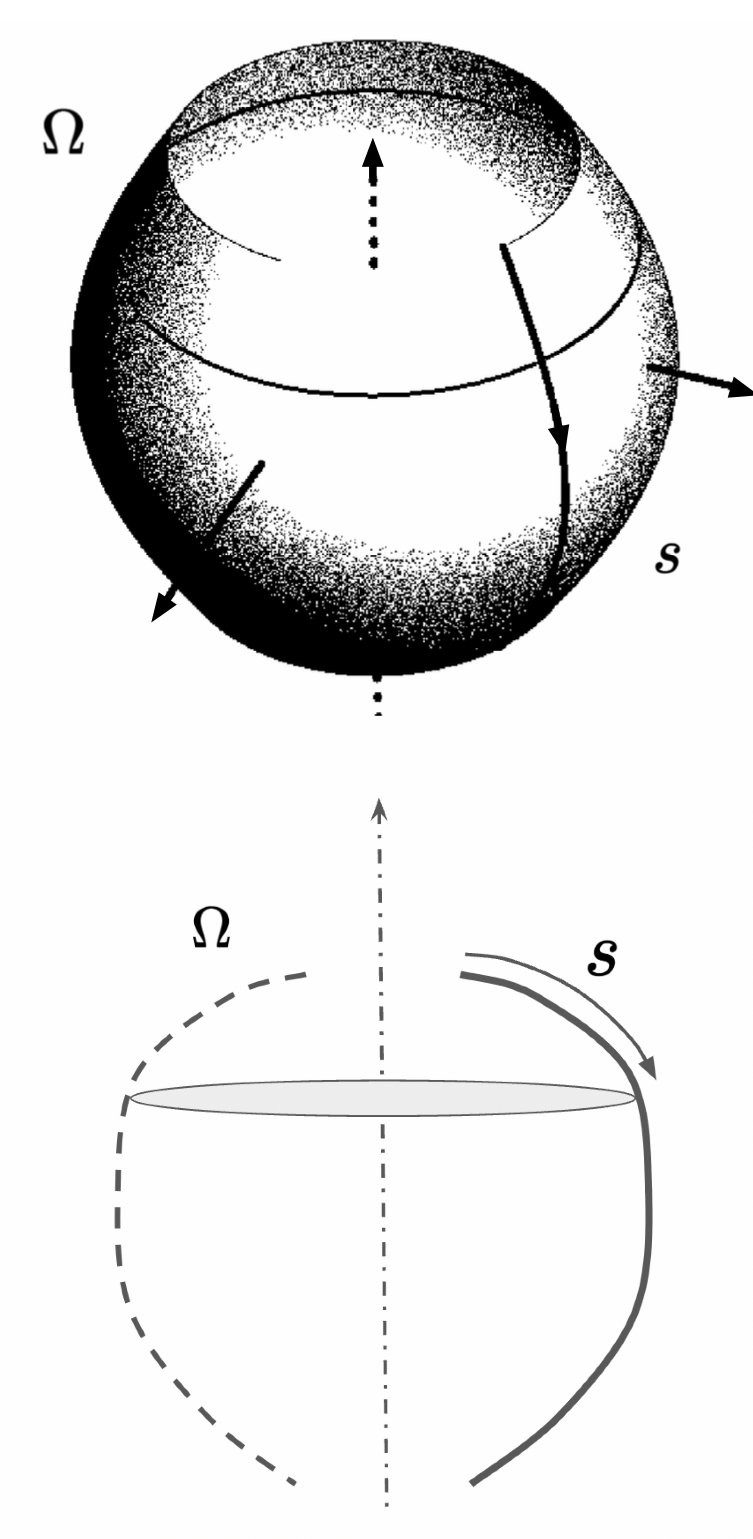} & 
\hspace{.1\textwidth} \includegraphics[width=.5\textwidth]{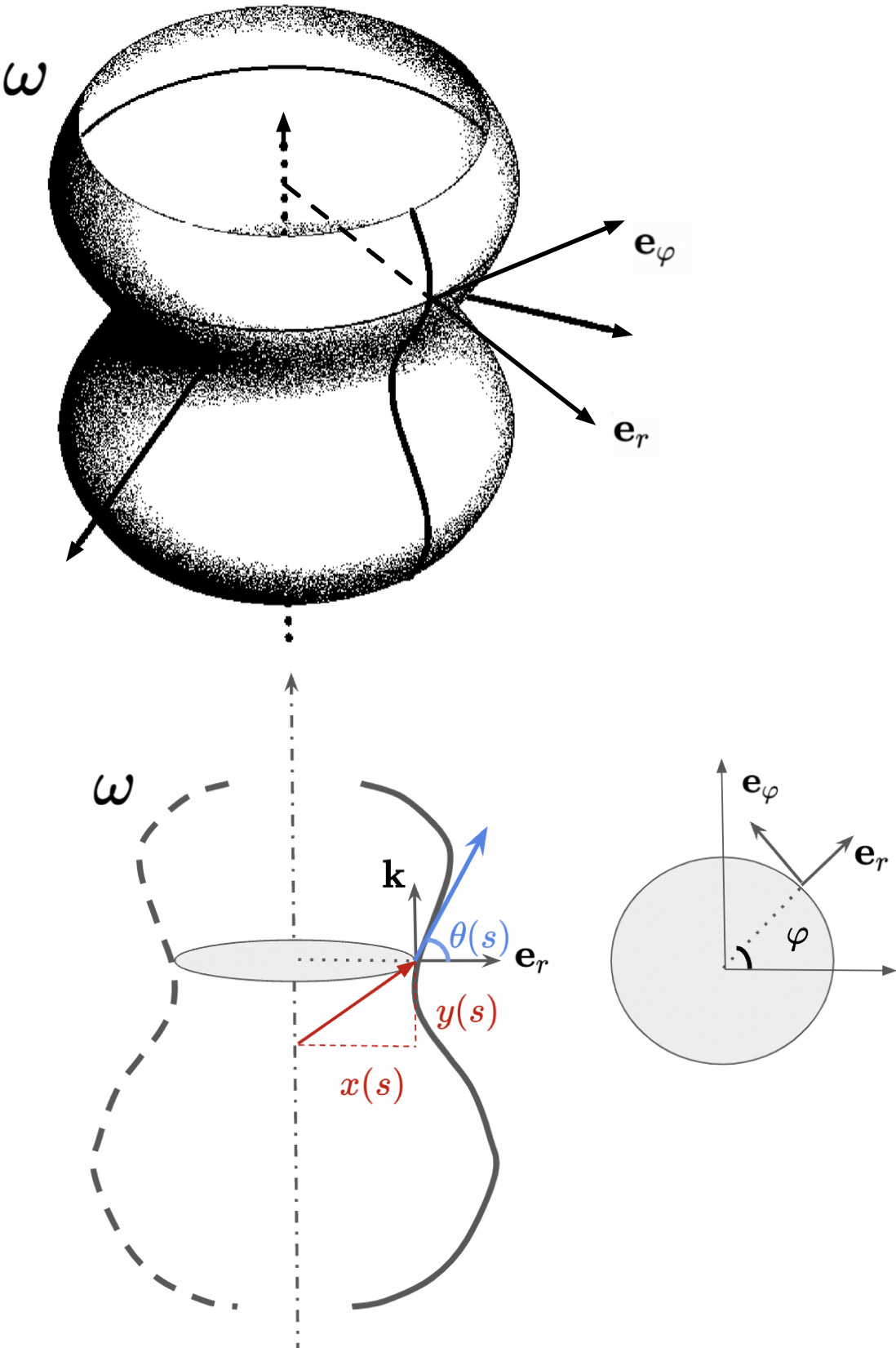}\\
(a) & (b)
\end{tabular}
\caption{(a) Reference configuration for an axisymmetric surface $\Omega$ (b) Deformed configuration $\omega$, with the deformation map $\mathbf{f}(s,\varphi)=x(s)\ee_r+y(s)\mathbf{k}$.}
\label{fig:schematic axi}
\end{figure}
Since this paper focuses on one-dimensional domains, we specialize {both above} theories, i.e., involving \eqref{mintheory} and \eqref{HCenergy}, for the case of axisymmetric configurations. Let $\Omega$ represent an axisymmetric reference configuration of the surface whose associated arc length parameter is denoted as $s$ (shown schematically in Fig.~\ref{fig:schematic axi}a). The total arc length along the cross-section is $L$, thus $s\in[0,L]$. Let $\omega$ represent the deformed configuration and the deformation map (shown {schematically} in red in Fig.~\ref{fig:schematic axi}b) for the surface $\omega$ is parametrized as 
\begin{equation}
\mathbf{f}(s,\varphi)=x(s)\ee_r(\varphi) + y(s)\mathbf{k},\;s\in[0,L],
\label{eq:axisymmetric para}
\end{equation}
where $\ee_r(\varphi)$, and $\mathbf{k}$ are the basis vectors in the cylindrical coordinates, with $\varphi$ denoting the azimuthal angle. The angle between the tangent to $\omega$ at $s$ and $\ee_r$ is denoted as $\theta(s)$. It is clear that
\begin{equation}
\tan\theta(s) = \frac{y'(s) }{x'(s) },
\notag
\end{equation}
where {the} prime denotes the derivative with respect to $s$. It can be shown 
that ({to avoid clutter,} omitting to write the explicit dependence on $s$)
\begin{equation}
da = x\sqrt{x'^2+y'^2}\;ds\;d\varphi,
\notag
\end{equation}
\[
H = \frac{1}{2}\Big(\frac{\theta'}{\sqrt{x'^2+y'^2}}+\frac{\sin\theta}{x}\Big),
\quad K = \frac{\theta'\sin\theta}{x \sqrt{x'^2+y'^2}}.
\]
{With restriction to consider only} axisymmetric configurations, the minimal surface energy \eqref{mintheory} can be written as 
\begin{equation}
\mathcal{E}_{MS} = 2\pi\int_0^L x\sqrt{x'^2+y'^2}\;ds,
\label{eq:axi MS}
\end{equation}
{while} the Helfrich-Canhman energy \eqref{HCenergy} can be written as
\begin{equation}
\begin{split}
\mathcal{E}_{HC} &= 2\pi\int_0^L \Bigg(\frac{\kappa}{4}\Big(\frac{\theta'}{\sqrt{x'^2+y'^2}}+\frac{\sin\theta}{x}\Big)^2\\
&\qquad + \kappa_g \frac{\theta'}{\sqrt{x'^2+y'^2}}\frac{\sin\theta}{x}\Bigg)x\sqrt{x'^2+y'^2}\;ds.
\label{eq:axi HC}
\end{split}
\end{equation}
{By inspection} it is clear that \eqref{eq:axi MS} is a first-order theory since the highest derivative appearing in the functional is of first-order, while \eqref{eq:axi HC} is a second-order theory {where} the functional depends on $\theta' = (x'y''-y'x'')/(x'^2+y'^2)$ (see \ref{appx:sec:axisymm formulation} {of supplementary for details}). 

{
As another example of second-order theory, which is inherently one-dimensional, we consider the Euler Elastica \cite{euler1952methodus}.
In this classic example (chapter 19, \cite{love1927treatise}), a rod (of reference length $L$) is described as a curve in two dimensions with its deformation map parametrized as 
\begin{equation}
\mathbf{r}(s):=x(s)\mathbf{i}+y(s)\mathbf{j},
\label{eq:defrsxsys}
\end{equation}
where $s$ is a parametrization of the reference arc-length of the rod. The energy density (in the deformed configuration) for an Euler elastica is proportional to the square of the curvature, the latter explicitly given by $\kappa(s) = (x'y''-y'x'')/(x'^2+y'^2)^{3/2}$. Thus the total energy is given by
\begin{equation}
\mathcal{E}_{e} = \int_0^L \Big(\frac{x'y''-y'x''}{(x'^2+y'^2)^{3/2}}\Big)^2\;\sqrt{x'^2+y'^2}\;ds.
\label{eq:euler}
\end{equation}
Indeed, \eqref{eq:euler} belongs to the class of second-order theories.
}

As mentioned before, {these three} examples serve as templates for the discussion below {where we follow a pedestrian approach, mostly employing freshman calculus, in the classical framework of the calculus of varitions \cite{Gelfand} and restrict attention to functionals on one-dimensional domains}.
We remark that, in addition to the surface of revolution described above, the extruded surfaces, obtained from a curve, can also be considered (where $K=0$) leading to another class of simple one dimensional problems but we omit the details for these as they appear to be less interesting from viewpoint of applications we have discussed so far.

\section*{Notation}
{Following the standard notation in mathematics and mechanics \cite{gurtin1982introduction},} the set of real numbers is denoted by $\mathbb{R}$ whereas the $n$ dimensional real linear space (with $n$-tuplets) is denoted by $\mathbb{R}^n$, which is further equipped with the standard inner-product {which we denote with a dot, i.e., $\mathbf{a}\cdot\mathbf{b}$}. 
{The standard basis of $\mathbb{R}^n$ is denoted by $\{\ee_i\}_{i=1}^n.$}
We use boldface font to represent the vectors and tensors. 
{The inherited standard basis for second order tensors is $\{\ee_i\otimes\ee_j\}_{i,j=1}^n$ with the second order tensor $\mathbf{a}\otimes\mathbf{b}$, for given $\mathbf{a},\mathbf{b}\in\mathbb{R}^n$, defined as a linear transformation from $\mathbb{R}^n$ to $\mathbb{R}^n$ such that $(\mathbf{a}\otimes\mathbf{b})[\mathbf{x}]=\mathbf{a}(\mathbf{x}\cdot\mathbf{b})$ for all $\mathbf{x}\in\mathbb{R}^n$; the square brackets denote the linear action of an operator, which is mostly clear from the context.}
We also employ the usual summation convention when the indicial notation is unambiguous to describe tensorial manipulations, {for example $\mathbf{a}=a_i\ee_i, \mathbf{A}[\mathbf{x}]=\mathbf{A}\mathbf{x}=A_{ij}x_j\ee_i, \mathbf{a}\cdot\mathbf{b}=a_ib_i,$ etc for vectors $\mathbf{a}, \mathbf{b}, \mathbf{x}$ and second order tensor $\mathbf{A}$}; otherwise, we stick to explicit sums. 
{In particular $\mathbf{a}\otimes\mathbf{b}=a_ib_j\ee_i\otimes\ee_j.$ We define the transpose of a tensor in a way that $(\mathbf{a}\otimes\mathbf{b})^T=a_jb_i\ee_i\otimes\ee_j.$}
$D_s$ represents the total derivative with respect to $s$ treating all arguments as functions of $s$. $D_1, D_2, \dotsc$ represents the (partial) derivative with respect to $1$st argument, $2$nd argument. For brevity and emphasis, we sometimes abuse the notation and alternate between 
{a comma notation and an explicit form, i.e., $
f_{,x}(x,y)=\partial_x f(x,y)=\frac{\partial}{\partial x} f(x,y)$, $
f_{,xy}= \partial_{xy}f(x,y)=\frac{\partial^2}{\partial y\partial x} f(x,y)$, etc. }
We use square brackets to represent a linear operator's action with an exception of the similar notation for a closed interval of the real line.
The unit sphere in $n$-dimensions is $\mathcal{S}^{n-1},$
while the corresponding torus $\mathcal{S}^1\times\dotsc\times\mathcal{S}^1$ is denoted by $(\mathcal{S}^1)^{n-1}$.

In the rest of this paper, we only consider functionals on one-dimensional domains, which we set to $[0,1]$ without any loss of generality. We assume that the Lagrangian $\mathscr{L}$ is, in general, a smooth function of $s\in[0,1]$, an $n$-dimensional field
\begin{equation}
{\uu(s) :=u_i(s)\ee_i\in\mathbb{R}^n,}
\label{eq:defuui}
\end{equation}
{i.e., simply the $n$-tuplet $(u_1(s),u_2(s),\dotsc,u_n(s))$,} and its first and second derivatives, $\uu{'}(s)$ and $\uu''(s)$.
The energy functional {with Lagrangian $\mathscr{L}$} is given by
\begin{equation}
\mathcal{E}(\uu)=\int_0^1\mathscr{L}(s,\uu(s),\uu'(s),\uu''(s))\;ds.
\label{eq:not:functional and lagrangian}
\end{equation}
If $\mathscr{L}$ explicitly depends on $\uu''$, we call it a \emph{second-order Lagrangian}. We assume that $n>1$ and denote the Euler-Lagrange operator acting on $\uu$ as
\begin{subequations}
\begin{equation}
\mathfrak{E}_\mathscr{L}(\uu):=\mathscr{L}_\uu - D_s\mathscr{L}_{\uu'}+D_{ss}\mathscr{L}_{\uu''},
\label{def:eulerop}
\end{equation}
where {
\begin{equation}
\begin{split}
\mathscr{L}_{\uu}&:=\mathscr{L}_{,z_i}(s,\mathbf{z},\mathbf{p},\mathbf{q})|_{\mathbf{z}=\uu(s),\mathbf{p}=\uu'(s),\mathbf{q}=\uu''(s)}\ee_i,\\
\mathscr{L}_{\uu'}&:=\mathscr{L}_{,p_i}(s,\mathbf{z},\mathbf{p},\mathbf{q})|_{\mathbf{z}=\uu(s),\mathbf{p}=\uu'(s),\mathbf{q}=\uu''(s)}\ee_i,\\
\mathscr{L}_{\uu''}&:=\mathscr{L}_{,q_i}(s,\mathbf{z},\mathbf{p},\mathbf{q})|_{\mathbf{z}=\uu(s),\mathbf{p}=\uu'(s),\mathbf{q}=\uu''(s)}\ee_i.
\end{split}
\label{def:Luetc}
\end{equation}
} 
\end{subequations}
The first variation \cite{Gelfand} of \eqref{eq:not:functional and lagrangian} with respect to smooth perturbations $\delta\uu$ of $\uu$ is denoted as $\delta\mathcal{E}[\delta\uu]$, {a linear functional of $\delta\uu$}. That is,
\begin{equation}
\delta\mathcal{E}[\delta\uu] = \int_0^1 (\mathscr{L}_{\uu}\cdot\delta\uu + \mathscr{L}_{\uu'}\cdot\delta\uu' + \mathscr{L}_{\uu''}\cdot\delta\uu'') \;ds,
\label{eq:not:first variation}
\end{equation}
where the dot represents the Euclidean inner product in $\mathbb{R}^n$.

{Let us introduce} the notation 
\begin{equation}
\mathcal{C}^m :=C^m([0,1];\mathbb{R}^n
\label{eq:Cmdef}
\end{equation}
We assume that $\uu\in\mathcal{C}^4$ for any non-negative integer $m$) and the perturbations
\begin{equation}
\delta\uu\in\mathcal{C}^\infty_c:=C^\infty_c([0,1];\mathbb{R}^n). 
\label{eq:deluCc}
\end{equation}
As is conventional in geometry and analysis, we use the word `curve' to describe a one-dimensional manifold in $\mathbb{R}^{n}$ where $n$ is a positive integer; we use $n=2$ for illustration of some results. The function space
\begin{align}
\mathcal{C}^m_R:=\{\uu\in \mathcal{C}^m: \text{the curve } \uu([0,1])\subset\mathbb{R}^{n}\text{ admits regular parametrization}\}
\label{def:C2R}
\end{align}
is used frequently in the paper. 
Notice that $\mathcal{C}^m_R$ is not a subspace of $\mathcal{C}^m$ but it can be considered as a candidate function space for variational analysis as $\mathcal{C}_c^\infty$ variations do not alter regularity mentioned in \eqref{def:C2R}. 

{In the context of the statements and discussion in the following, specially using \eqref{eq:not:functional and lagrangian} and \eqref{eq:not:first variation}, we emphasize that whenever} $\mathscr{L}$ is independent of $\uu''$, we call such it a \emph{first-order Lagrangian}; {Naturally,} in this case, it is sufficient to assume that $\uu\in \mathcal{C}^2$.

\section{Revisiting {parametrization} invariance}
\label{sec_revisit}
The purpose of this section is to recall an important {\em consequence} of {parametrization} invariance. The reader may note that we do not claim any originality in this section {but provide these details as a simple recall of certain elementary facts}.

We consider a diffeomorphism $\varpi:[0,1]\to[0,1]$ 
satisfying $\varpi(0)=0$ and $\varpi(1)=1$.
As a \emph{parametrization of the curve defined by $\uu$}, we use $t=\varpi^{-1}(s)$ as the symbol for a new parameter in place of $s$ for given $\uu\in\mathcal{C}^4_R$; thus, $s=\varpi(t)$ such that $\varpi(0)=0, \varpi(1)=1$. 
\begin{definition}[{Parametrization} Invariance]
$\mathcal{E}$ is said to be \emph{{parametrization} invariant} if for any $\uu\in\mathcal{C}_R^4$
\begin{align}
\mathcal{E}(\uu\circ\varpi)=\mathcal{E}(\uu),
\label{eq:parmiv1D}
\end{align}
for any parametrization $\varpi$.
\label{defparinvar}
\end{definition}

As a consequence of this symmetry, the tangential component of the Euler-Lagrange operator vanishes trivially. This is stated formally by 
\begin{prop}
If the functional \eqref{eq:not:functional and lagrangian} is {parametrization} invariant, then for any $\uu\in \mathcal{C}^4_R$,
\begin{equation}
\mathfrak{E}_{\mathscr{L}}(\uu)\cdot\uu'\equiv 0,
\label{eq:para inv1 EL}
\end{equation}
{that is, $(\mathscr{L}_\uu- D_s\mathscr{L}_{\uu'}+D_{ss}\mathscr{L}_{\uu''})\cdot\mathbf{u}'\equiv 0$.} 
\label{lemma2_2}
\end{prop}
This result is a special case of Noether's second theorem \cite{olver1986noether}{, and a proof of this result is included in the supplementary material.} For first-order Lagrangians, in the above proposition, naturally it sufficient to assume that $\uu\in \mathcal{C}^2_R$.

Indeed, the (integral or global) condition of {parametrization} invariance {constrains} a first-order Lagrangian to take a special form. 
{For the convenience of the reader, we} present a classical{, and an elementary,} argument \cite{GiaquintaV1} for the more general $n$-dimensional case.
Due to \eqref{eq:parmiv1D}, we have the equations:
\begin{align}
\int_{s_1}^{s_2} {\mathscr{L}}(s, \uu(s),\uu'(s))\;ds&=\int_{\varpi^{-1}(s_1)}^{\varpi^{-1}(s_2)} {\mathscr{L}}\Big(\varpi(t), \uu\circ\varpi(t),(\uu\circ\varpi){'}(t)/\varpi'(t)\Big)\varpi'(t)\;dt\notag\\
&=\int_{\varpi^{-1}(s_1)}^{\varpi^{-1}(s_2)} {\mathscr{L}}\Big(t, \uu\circ\varpi(t),(\uu\circ\varpi)'(t))\Big)\;dt,
\label{eqL1}
\end{align}
{with $(s_1,s_2)\subset[0,1].$}
The first {equality in \eqref{eqL1} arises due to} a change of variable of integration, while the second equation is a statement of invariance.
Using {the arbitrariness of $s_1,s_2$, the condition equivalent to {parametrization} invariance is found to be}
\[
{\mathscr{L}}\Big(\varpi(t), \uu\circ\varpi(t),(\uu\circ\varpi){'}(t)/\varpi'(t)\Big)\varpi'(t)
={\mathscr{L}}(t, \uu\circ\varpi(t),(\uu\circ\varpi)'(t)))
\]
{for each $t\in(0,1)$.} It is {well known} \cite{GiaquintaV1} that the following two conditions:
(1) ${\mathscr{L}}$ is independent of {$s$} and 
(2) ${\mathscr{L}}$ is a homogeneous function of degree one in $\uu'(t)$, 
are necessary and sufficient for the integral form of {parametrization} invariance as stated in Definition {\ref{defparinvar}}.
We {formalize} this observation as
\begin{theorem}
A functional $\mathcal{E}$ with a first-order Lagrangian is {parametrization} invariant according to Definition {\ref{defparinvar}}
if and only if the Lagrangian ${\mathscr{L}}$ is independent of the parameter $s$ and ${\mathscr{L}}$ is a homogeneous function of degree one in $\uu'$.
\label{thm:para inv lagrangian nd}
\end{theorem}

{As a counterpart of above statement for a second-order Lagrangian \cite{bates2016elastica} (for example, see Theorem A.31 and Remark A.32), we have
\begin{theorem}
A functional $\mathcal{E}$ with a second-order Lagrangian is {parametrization} invariant
if and only if
\begin{subequations}
\begin{eqnarray}
\uu'\cdot\mathscr{L}_{\uu''}&=&0, \label{para2a} \\
\uu'\cdot(\mathscr{L}_{\uu'}-\frac{d}{ds}\mathscr{L}_{\uu''}) +\uu''\cdot\mathscr{L}_{\uu''}-\mathscr{L}&=&0. \label{para2b}
\end{eqnarray} 
\end{subequations}
\label{thm:2nd:para inv lagrangian nd}
\end{theorem}
}
\begin{figure}
\centering
\includegraphics[width=.6\linewidth]{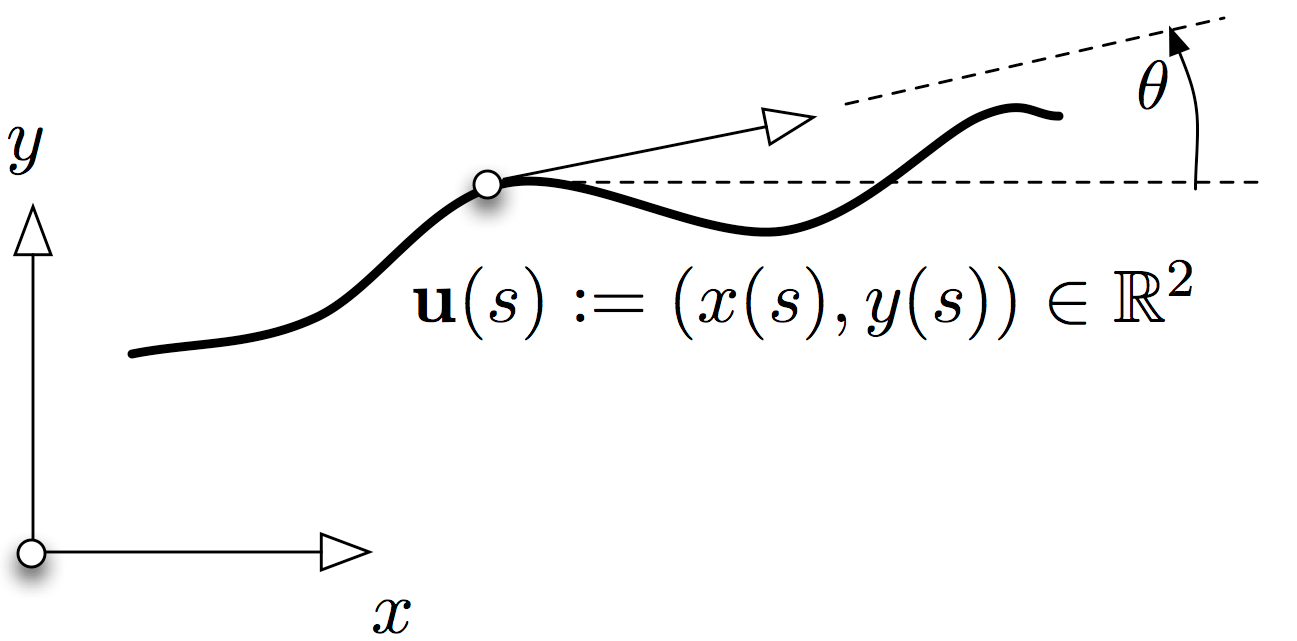}
\caption{Schematic of the angle $\theta$ in $\mathbb{R}^2$.}
\label{fig:schematic1}
\end{figure}

\begin{remark}
In the case of $n=2$, i.e., when $\uu([0,1])=\{(x(s),y(s)): s\in[0,1]\}\subset\mathbb{R}^2$ it can be shown {for first order case} that \cite{GiaquintaV1}
\begin{equation}
\mathscr{L}(s,(x(s),y(s)),(x'(s),y'(s)))=f\Big(x(s),y(s),\theta(s)\Big)\sqrt{x'(s)^2+y'(s)^2}
\label{eq:para inv lagrangian 2d}
\end{equation}
for some smooth function $f$, where $\theta(s)=\arctan(y'(s)/x'(s))$ can be interpreted as the tangent angle to the curve and $s\in[0,1]$; see Fig.~\ref{fig:schematic1} (c.f., Fig.~\ref{fig:schematic axi}b). This result is a consequence of Euler's theorem for homogeneous functions. 
{Similarly,} for the two-dimensional (i.e., $\uu(s)=(x(s),y(s))$) problem with a second-order Lagrangian \cite{GiaquintaV1}, {the Theorem \ref{thm:2nd:para inv lagrangian nd} translates into the statement:}
A functional $\mathcal{E}$ with a second-order Lagrangian is {parametrization} invariant
if and only if the Lagrangian has the form:
\begin{equation}
\begin{split}
&\mathscr{L}
{(s,(x(s),y(s)),(x'(s),y'(s)),(x''(s),y''(s)))}\\
&= g\Big(x(s),y(s),\theta(s), \theta'(s)/\sqrt{x'(s)^2+y'(s)^2}\Big)\sqrt{x'(s)^2+y'(s)^2},
\label{eq:rep form 2nd order}
\end{split}
\end{equation}
where $g$ is an arbitrary (smooth) function and $\theta(s)=\arctan(y'/x')$.
{We retain the presence of $s$ on the left hand side of \eqref{eq:para inv lagrangian 2d} and \eqref{eq:rep form 2nd order} for emphasis.}
At this point, it is appropriate to compare the structure of 
{\eqref{eq:para inv lagrangian 2d} with that of \eqref{eq:axi MS},
and \eqref{eq:rep form 2nd order} with that of \eqref{eq:axi HC}. }
\end{remark}

In summary, we see in Proposition \ref{lemma2_2} that the natural {\em integral} form of {parametrization} invariance implies {a form of {\em differential}} condition that the tangential component of the Euler-Lagrange operator vanishes. {Parametrization} invariance also restricts the form of the Lagrangian (c.f., Theorems \ref{thm:para inv lagrangian nd} and \ref{thm:2nd:para inv lagrangian nd}). In this work, we explore a connection, {in a sort of opposite implication}, between 
Proposition \ref{lemma2_2}
and Theorems \ref{thm:para inv lagrangian nd} and \ref{thm:2nd:para inv lagrangian nd}, {i.e.} without the hypothesis of {parametrization} invariance.

\section{\texorpdfstring{{$\mathcal{T}$-}}{T} Lagrangian: case of tangential variations for first-order Lagrangians}
\label{sec_tangentvar}
The concept of {$\mathcal{T}$-}{}Lagrangian that we discuss below is motivated by relaxing the notion of {parametrization} invariance. Instead of requiring a \emph{global} (integral) criterion, the invariance condition \eqref{eq:parmiv1D}, we enforce a weaker, \emph{local} (differential) condition of {$\mathcal{T}$-}{}Lagrangian, that we precisely state below.

If $\uu\in \mathcal{C}^2_R$ {(recall \eqref{def:C2R})}, we define variations that are tangential and normal to the curve;
note that $\uu'$ is non-vanishing.
That is,
\begin{definition} [Normal and Tangential variations]
Assume $n>1.$
Given a $\uu\in\mathcal{C}_R^2$ the set
\begin{equation}
N_\uu :=\{ \mathbf{v} \in \mathcal{C}^\infty_c: 
{\mathbf{v}\cdot \mathbf{u'}=0}\}
\notag
\end{equation}
is the set of \emph{normal variations} and 
\begin{equation}
T_\uu :=N_\uu^\perp=\{ \mathbf{v} \in \mathcal{C}^\infty_c: 
{\mathbf{v}=\psi\mathbf{u'}, \text{for some }\psi\in\mathcal{C}^\infty_c}
\},
\notag
\end{equation}
is the set of \emph{tangential variations}. 
\end{definition}
Clearly, {the direct sum} $T_\uu\oplus N_\uu$ coincides with $\mathcal{C}^\infty_c$.
Thus, $T_\uu$ contains functions that have components tangential to $\uu$, and $N_\uu$ contains functions with components normal to $\uu$.

Let us consider the following condition, motivated by \eqref{eq:para inv1 EL}, stated as
\begin{definition}[{$\mathcal{T}$-}{}Lagrangian]
${\mathscr{L}}$ is a \emph{{$\mathcal{T}$-}{}Lagrangian} if for every $\uu\in \mathcal{C}_R^2$, $\delta\mathcal{E}[\delta\uu]=0$ for all $\delta\uu\in T_\uu$.
\label{def:T lag}
\end{definition}
That is, functional with a {$\mathcal{T}$-}{}Lagrangian {necessitates that} its first variation {is} trivially zero with respect to variations in a direction tangential to {the curve} $\uu$. This definition contrasts with that of a classical null Lagrangian \cite{Rund66} whose first variation is trivially zero for \emph{all} variations $\delta\uu(s)$. Thus, {$\mathcal{T}$-}{}Lagrangian can be interpreted as a weakened notion of the classical null Lagrangian; recall that the Euler-Lagrange operator for a null Lagrangian is trivial. We show that an equivalent condition to check for {the defining property of} a {$\mathcal{T}$-}{}Lagrangian is a natural condition that the tangential component of the Euler-Lagrange operator is trivially zero.
\begin{prop}
${\mathscr{L}}(s,\uu,\uu')$ is a {$\mathcal{T}$-}{}Lagrangian if and only if the tangential component of the Euler-Lagrange operator vanishes for all $\uu\in \mathcal{C}^2_R$. That is,
\begin{equation}
\uu'\cdot\mathfrak{E}_\mathscr{L}(\uu)=\uu'\cdot(-D_s{{\mathscr{L}}}_{\uu'} + {{\mathscr{L}}}_\uu)\equiv 0.
\label{eq:tgtEL=0}
\end{equation}
\label{lem:pnlagrangian_lem}
\end{prop}
\begin{proof}
\begin{subequations}
We first prove the forward implication, i.e., if $\mathscr{L}$ is a {$\mathcal{T}$-}{}Lagrangian then \eqref{eq:tgtEL=0} holds. For any given $\uu\in\mathcal{C}^2_R$, let us consider variations $\delta\uu\in T_\uu$. The first variation is given by 
\begin{equation}
\delta \mathcal{E}[\delta\uu] = \int_0^1 \big({{\mathscr{L}}}_{\uu}\cdot \delta\uu + {{\mathscr{L}}}_{\uu'}\cdot\delta\uu'\big) ds.
\label{eq:deltaE}
\end{equation}
which after integrating by parts yields
\begin{equation}
\delta\mathcal{E}[\delta\uu]=\int_0^1\Big(\mathscr{L}_\uu-D_s\mathscr{L}_{\uu'}\Big)\cdot\delta\uu \;ds=0.
\label{eq:delta E=0}
\end{equation}
Since $\delta\uu\in T_\uu$ it has the form {$\delta\uu = \psi(s)\uu'$}, where $\psi\in C_c^\infty[0,1]$. Therefore,
\begin{equation}
\delta\mathcal{E}[\delta\uu]=\int_0^1\Big(\mathscr{L}_\uu-D_s\mathscr{L}_{\uu'}\Big)\cdot{\psi\uu'}\;ds=0.
\notag
\end{equation}

Since $\psi(s)$ is arbitrary, we obtain
\begin{equation}
\uu'\cdot \Big(\mathscr{L}_\uu-D_s\mathscr{L}_{\uu'}\Big) = 0.
\end{equation}
The converse follows from the fact that if \eqref{eq:tgtEL=0} holds, then the integrand in \eqref{eq:delta E=0} must be zero, whence $\delta\mathcal{E}[\delta\uu]=0$ for all admissible variations.
\end{subequations}
\end{proof}

Just as in the case of {parametrization} invariant functionals, the tangential component of the Euler-Lagrange operator of a {$\mathcal{T}$-}{}Lagrangian is, {thus}, trivially zero. But for the latter, the converse also holds. {In short}, the vanishing of the tangential component is an equivalent \emph{local} criterion for a {Lagrangian to be a} {$\mathcal{T}$-}{}Lagrangian.

\subsection{Representation theorem for first-order \texorpdfstring{{$\mathcal{T}$-}}{T} Lagrangian}
Recall that Theorem \ref{thm:para inv lagrangian nd} and Theorem \ref{thm:2nd:para inv lagrangian nd} characterize the form of {parametrization} invariant Lagrangians. In this section, we derive the representation theorem for {$\mathcal{T}$-}{}Lagrangians. We focus on first-order Lagrangians in this section, {while} we consider the second-order case {in the next}.

{In the proof of the representation theorem, a crucial role is played by}
\begin{lemma}
If ${\mathscr{L}}(s,\uu,\uu')$ is a {$\mathcal{T}$-}{}Lagrangian then ${\mathscr{L}}$ has the form
\begin{subequations}
\begin{equation}
{\mathscr{L}}(s,\uu,\uu') = \uu' \cdot {\mathbf{A}}(s,\uu,{\boldsymbol{\theta}})+ C(s,\uu),
\label{eq:rep lemma eq ndof}
\end{equation}
where $C$ is an arbitrary function (of variables indicated), {${\mathbf{A}} := 
A_i\ee_i$}, with $A_i$ $(i=1,\dotsc,n)$ as functions of $s, \uu$, and ${\boldsymbol{\theta}}\in (\mathcal{S}^1)^{n-1}$, {such that}
\begin{equation}
{\mathbf{A}}^T_{{\boldsymbol{\theta}}}(s,\uu,{\boldsymbol{\theta}})[\uu']=\mathbf{0};\text{ for all }\uu\in \mathcal{C}^2_R,
\label{eq:pA_t=-qB_t ndof}
\end{equation}
where {$\mathbf{A}_{{\boldsymbol{\theta}}}=\frac{\partial A_i}{\partial\theta_j} \ee_i\otimes\ee_j$,} and 
\label{lem:L form ndof}
${\boldsymbol{\theta}}=(\theta_1, \theta_2, \dotsc, \theta_{n-1})$ is associated with $\hat{\p}=\uu'/\norm{\uu'}\in\mathcal{S}^{n-1}$ in the following manner (schematically shown in Fig.~\ref{Fig2}): 
\begin{equation}
\begin{array}{l} 
\hat{p}_1=\cos \theta_1,\\
\hat{p}_2=\sin \theta_1\cos \theta_2,\\
\hat{p}_3=\sin \theta_1\sin \theta_2\cos \theta_3,\\
\vdots\\
\hat{p}_{n-1}=\sin \theta_1\sin \theta_2\sin \theta_3\dotsc\sin\theta_{n-2}\cos\theta_{n-1},\\
\hat{p}_n=\sin \theta_1\sin \theta_2\sin \theta_3\dotsc\sin\theta_{n-2}\sin\theta_{n-1}.
\end{array}
\label{eq:pn equations}
\end{equation}

\begin{figure}
\centering
\includegraphics[width=.45\linewidth]{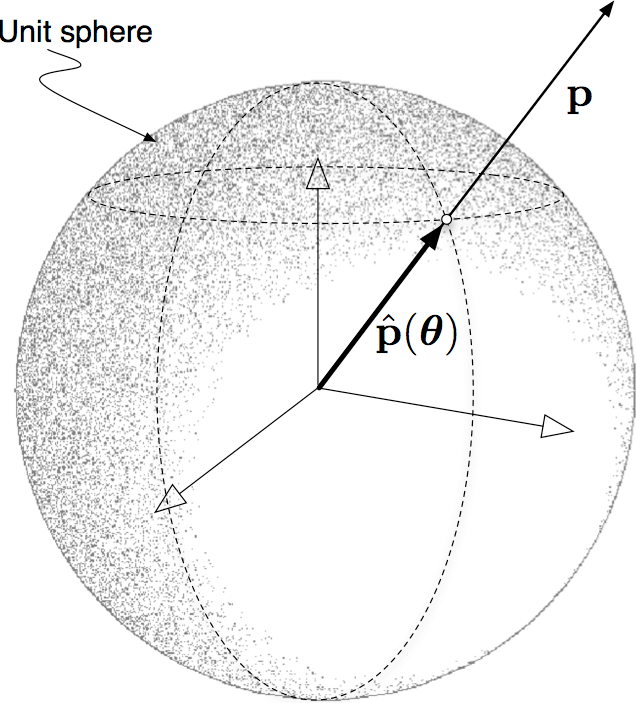}
\caption{Schematic of ${\hat{\p}}({\boldsymbol{\theta}})$ in the polar coordinates $\p=\p(r,{\boldsymbol{\theta}}) = r{\hat{\p}}({\boldsymbol{\theta}})$ in $\mathbb{R}^n$.}
\label{Fig2}
\end{figure}

\end{subequations}
\end{lemma}
\begin{proof}
\begin{subequations}
Since ${\mathscr{L}}(s,\uu,\uu')$ is a {$\mathcal{T}$-} Lagrangian, we have (from Proposition
\ref{lem:pnlagrangian_lem}) for any $\uu\in \mathcal{C}^2_{R}$,
\begin{equation*}
(-D_s({{\mathscr{L}}}_{\uu'}) + {{\mathscr{L}}}_\uu)\cdot\uu'=0.
\end{equation*}
Let us define 
\[
\p:=\uu',
\]
i.e., we consider ${\mathscr{L}}(s,\uu,\p)$. 
{Recall the definitions (or consider their natural analogues) \eqref{def:Luetc}, i.e. 
\[
{\mathscr{L}}_{\p}:={\mathscr{L}}_{,p_i}\ee_i, \quad
{\mathscr{L}}_{\p\p}:={\mathscr{L}}_{,p_ip_j}\ee_i\otimes\ee_j, \quad
{\mathscr{L}}_{\p\uu}:={\mathscr{L}}_{,p_iu_j}\ee_i\otimes\ee_j.
\]}
Expanding $D_s{\mathscr{L}}_{\p}$ in the previous equation, we have
\begin{equation*}
(-{\mathscr{L}}_{\p\p}[\p']-{\mathscr{L}}_{\p\uu}[\p]-{\mathscr{L}}_{\p s} + {\mathscr{L}}_{\uu})\cdot\p= 0.
\end{equation*}
Collecting the $\p'$ terms, we get
\begin{equation*}
-{\mathscr{L}}_{\p\p}[\p']\cdot\p + (\text{ terms that depend on }s, \uu, \p)= 0.
\end{equation*}
Since this must hold for all $\uu\in\mathbb{R}^n$ and $\p\in\mathbb{R}^n$, the coefficients
of $\p'$ must vanish:
\begin{equation}
{\mathscr{L}}_{\p\p}[\p]=\mathbf{0}.
\label{eq:Lpp conditions ndof}
\end{equation}
By Euler's theorem for homogeneous functions, \eqref{eq:Lpp conditions ndof}
hold iff the components of ${\mathscr{L}}_{\p}$ are
homogeneous functions of $\p$ of degree zero. That is, for any $\lambda>0$,
\begin{equation}
{\mathscr{L}}_{\p}(s,\uu, \lambda\p)= {\mathscr{L}}_{\p}(s,\uu, \p).
\label{eqs:scale invariance ndof}
\end{equation}
It then follows (see Lemma \ref{lem:consequence of homogeneity} in the appendix) that
\begin{equation}
{\mathscr{L}}_{\p}(s,\uu, \p) = {\mathbf{A}}(s,\uu,{\boldsymbol{\theta}}),
\label{eq:AB ndof}
\end{equation}
where {$\mathbf{A} = A_i\ee_i$,} for some functions $A_i, i=1, \dotsc, n$ and ${\boldsymbol{\theta}}\in (\mathcal{S}^1)^{n-1}$. {Note that $\bm{\theta}$ is related to $\mathbf{p}$ through the relation 
\begin{equation}
\mathbf{p} =\mathbf{p}(r,\bm{\theta}) := r\hat{\mathbf{p}}(\bm{\theta})
\label{eq:pprptheta}
\end{equation}
as shown in Fig.~\ref{Fig2}, c.f., \eqref{eq:pn equations}.}

Differentiating {\eqref{eq:AB ndof} with respect to $\p$ and contracting the result with $\mathbf{p}$, we obtain, in indicial notation,}
\begin{equation}
p_j {\mathscr{L}}_{p_ip_j}=rA_{i,r}=0,
\label{eq:pL_pp=0}
\end{equation}
{where the first equation follows from \eqref{frftheta} (see appendix) and second from \eqref{eq:AB ndof}. Again from \eqref{frftheta}, we find that
\begin{equation}
A_{i,{\boldsymbol{\theta}}}=({\p}_{{\boldsymbol{\theta}}})^T[{\mathscr{L}}_{p_i,\p}],
\label{eq:Aibtheta}
\end{equation}
which in components is given by $A_{i,\theta_j}={p}_{k,\theta_j}({\mathscr{L}}_{p_ip_k})$. 
Thus, using \eqref{eq:pL_pp=0} and \eqref{eq:Aibtheta}, we deduce that $A_{i,\theta_j}p_i={p}_{k,\theta_j}({\mathscr{L}}_{p_ip_k})p_i=0$ (since ${\mathscr{L}}_{p_ip_k}={\mathscr{L}}_{p_kp_i}$). That is,}
\begin{equation}
{\mathbf{A}}_{,{\boldsymbol{\theta}}}^T(s,\uu,{\boldsymbol{\theta}})[\p]=\mathbf{0}.
\notag
\end{equation}

Let us expand $\frac{\partial}{\partial\lambda}{\mathscr{L}}(s,\uu,\lambda \p)$ as follows:
\begin{equation*}
\frac{\partial}{\partial\lambda}{\mathscr{L}}(s,\uu,\lambda \p) = \p\cdot {\mathscr{L}}_{\p}(s,\uu,\lambda \p).
\notag
\end{equation*}
Due to \eqref{eq:AB ndof}, we have
\begin{equation}
\frac{\partial}{\partial\lambda}{\mathscr{L}}(s,\uu,\lambda \p) = \p\cdot {\mathbf{A}}(s,\uu, {\boldsymbol{\theta}}).
\notag
\end{equation}
Integrating this equation with respect to to $\lambda$ { (from 0 to 1), we obtain
\begin{equation}
{\mathscr{L}}(s,\uu,\p) = \p\cdot {\mathbf{A}}(s,\uu, {\boldsymbol{\theta}})+ C(s,\uu).
\label{eq:formofL ndof}
\end{equation}
where $C(s,\uu):=\mathscr{L}(s,\uu,\mathbf{0})$} and ${\mathbf{A}}$ satisfies \eqref{eq:pA_t=-qB_t ndof}.
\end{subequations}
\end{proof}

To simplify this result further, let us first recall the special case
of Poincar\'e Lemma applied to $2n$ dimensions, i.e.,
\begin{lemma}[Poincar\'e Lemma]
A smooth 1-form on $\mathbb{R}^{2n}$,
\begin{equation*}
\alpha = c\;ds + \sum_{i=1}^n a_i du_i + \sum_{i=1}^{n-1} e_i d\theta_i,
\end{equation*}
where $a_i,e_i,c$ are functions of $(s,\uu,{\boldsymbol{\theta}})\in\mathbb{R}\times\mathbb{R}^{n}\times\mathbb{R}^{n-1}$, is closed ($d\alpha=0$) if and only if it is exact ($\alpha = d\psi$, for some $\psi(s,\uu,{\boldsymbol{\theta}})$). That is,
\begin{equation}
c_{,u_i} = a_{i,s},
c_{,\theta_i}=e_{i,s},
a_{i,u_j}=a_{j,u_i},
a_{i,\theta_j}=e_{j,u_i},
\iff c=\psi_{,s},\;a_i=\psi_{,u_i},e_i=\psi_{,\theta_i}.
\notag
\end{equation}
\label{lem:4d:poincare ndof}
\end{lemma}
As one of the main results of this paper, we state the representation theorem for a {$\mathcal{T}$-}{}Lagrangian as
\begin{theorem}
If ${\mathscr{L}}(s,\uu,\uu')$ is a {$\mathcal{T}$-} Lagrangian, then ${\mathscr{L}}$ has the local representation
\begin{equation}
{\mathscr{L}}\big(s,\uu(s),\uu'(s)\big)=D_s\Psi\big(s,\uu(s)\big) + \uu'(s)\cdot\mathbf{E}\big(\uu(s),{\boldsymbol{\theta}}(s)\big),
\notag
\end{equation}
where ${\boldsymbol{\theta}}\in (\mathcal{S}^1)^{n-1}$, {$\Psi$ is a (smooth) function of its arguments,} and $\mathbf{E}(\uu,{\boldsymbol{\theta}})$ is a vector of $n$ functions such that
\begin{equation}
E_{,{\boldsymbol{\theta}}}^T(\uu,{\boldsymbol{\theta}})[\uu']=\mathbf{0}.
\notag
\end{equation}
\label{thmatangt_full ndof}
\end{theorem}
\begin{proof}
In light of \eqref{eq:rep lemma eq ndof}, with $\p=\uu',$ we have
\begin{equation}
{\mathscr{L}}_{\p}{(s,\uu,\p)} = {\mathbf{A}}{(s,\uu,{\boldsymbol{\theta}})},\quad
{\mathscr{L}}_{\uu}{(s,\uu,\p)} = \p\cdot {\mathbf{A}}_{,\uu}{(s,\uu,{\boldsymbol{\theta}})} + C_{,\uu}{(s,\uu)}.
\notag
\end{equation}
Since ${\mathscr{L}}$ is a {$\mathcal{T}$-} Lagrangian, due to Proposition
\ref{lem:pnlagrangian_lem}, 
after simplifying \eqref{eq:tgtEL=0} further { and using \eqref{eq:pA_t=-qB_t ndof}}, we obtain
\begin{equation}
p_i(C_{,u_i}-A_{i,s})= 0.
\label{eq:tgt lag general temp ndof}
\end{equation}
Differentiating \eqref{eq:tgt lag general temp ndof} with respect to $p_i$ and using \eqref{eq:pA_t=-qB_t ndof}, we deduce that 
\begin{equation}
A_{i,s}{(s,\uu,{\boldsymbol{\theta}})}=C_{,u_i}{(s,\uu)}.
\label{eq:ABC ndof}
\end{equation}
Integrating these equations with respect to $s$, we obtain
\begin{equation}
A_i{(s,\uu,{\boldsymbol{\theta}})} = \int C_{,u_i}{(s,\uu)}\;ds + E_i(\uu,{\boldsymbol{\theta}}).
\label{eq:gen:a,b ndof}
\end{equation}
for some functions $E_i$ of $\uu$ and ${\boldsymbol{\theta}}$, and the {integral is an anti-derivative of $C_{,u_i}$} with respect to $s$, treating the other variables as constants. {Note that since any two anti-derivatives of $C_{,u_i}$ (with respect to $s$) differ by an arbitrary function of $\uu$ and ${\boldsymbol{\theta}}$, any one anti-derivative can be chosen in the above equation.} 

Note that since $C$ is independent of ${\boldsymbol{\theta}}$, we have $A_{i,{\boldsymbol{\theta}}}=E_{i,{\boldsymbol{\theta}}}$. Using \eqref{eq:pA_t=-qB_t ndof}, we obtain
\begin{equation}
E_{,{\boldsymbol{\theta}}}^T(\uu,{\boldsymbol{\theta}})[\p]=\mathbf{0}.
\label{eq:ct e_t = st f_t ndof}
\end{equation}
Let us define the following functions
\begin{equation}
a_i{(s,\uu,\p)}:=A_i{(s,\uu,\p)}-E_i{(s,\uu,\p)} = \int C_{,u_i}{(s,\uu)}\;ds,\;
e_i{(s,\uu,\p)}:=0,
\label{eqs:abcd ndof}
\end{equation}
where $0$ is the zero function (and the second equality in the first equation follows from \eqref{eq:gen:a,b ndof}). 

We verify that the conditions for Poincar\'e Lemma \ref{lem:4d:poincare ndof} are satisfied {for $C$, $a_i$, and $e_i$}:
{
\begin{subequations}
\begin{eqnarray}
a_{i,s} &=& C_{,u_i}\\
 C_{,\theta_i} &=& 0 = e_{i,s},\\
a_{i,u_j} &=& \int C_{,u_iu_j}\;ds = \int C_{,u_ju_i}\;ds = a_{j,u_i},\\
a_{i,\theta_j} &=& 0 = e_{j,u_i},
\end{eqnarray}
\end{subequations}
}
where the first two follow from the definition of $a_i$s (cf. \eqref{eqs:abcd ndof}), and the last two follow from the independence of $C$ with respect to ${\boldsymbol{\theta}}$ and the definition $e_i=0$.
It follows from Poincar\'e Lemma \ref{lem:4d:poincare ndof} that there exists a $\Psi(s,\uu,{\boldsymbol{\theta}})$ {such that}
\begin{equation}
{\mathbf{A}}-\mathbf{E} = \Psi_{,\uu},\; C =\Psi_{,s},\; \Psi_{,{\boldsymbol{\theta}}}=0.
\label{eqs:eqnsset1}
\end{equation}
In particular, due to the {third equation in above \eqref{eqs:eqnsset1}}, we see that $\Psi$ is independent of ${\boldsymbol{\theta}}$.

The general form for the {$\mathcal{T}$-} Lagrangian is thus given by
\begin{equation}
{\mathscr{L}}(s,\uu,\uu') = \uu'\cdot(\Psi_{,\uu}+\mathbf{E})+ \Psi_{,s},
\notag
\end{equation}
which can be simplified to 
\begin{equation}
{\mathscr{L}}\left(s,\uu(s),\uu'(s)\right)=\frac{d}{ds}\Psi\left(s,\uu(s)\right) + \uu'(s)\cdot\mathbf{E}\left(\uu(s),{\boldsymbol{\theta}}(s)\right),
\notag
\end{equation}
such that ${\boldsymbol{\theta}}$ and $\mathbf{E}$ satisfy the condition \eqref{eq:ct e_t = st f_t ndof}.
\end{proof}

\begin{remark}
When we specialize above theorem to the case $n=2$, i.e., with $\uu(s) = (x(s),y(s))$, we obtain the following form for the {$\mathcal{T}$-}{}Lagrangian $\mathscr{L}(s,x,y,x',y')$:
\begin{equation}
{\mathscr{L}}(s,x,y,x',y')=D_s\Psi(s,x(s),y(s)) + x'E_1(x,y,\theta) + y'E_2(x,y,\theta),\;
\label{eq:general form P2d}
\end{equation}
where $E_1(x,y,\theta)$ and $E_2(x,y,\theta)$ are functions such that
\begin{equation}
\cos\theta \frac{\partial E{_1}}{\partial\theta}(x,y,\theta) = \sin\theta \frac{\partial E_2}{\partial\theta}(x,y,\theta),
\notag
\end{equation}
and $\cos\theta=x'/\sqrt{x'^2+y'^2}$ and $\sin\theta=y'/\sqrt{x'^2+y'^2}$.
Observe that by multiplying and dividing the last two terms of \eqref{eq:general form P2d} by $\sqrt{x'^2+y'^2}$, these terms can be expressed as $\sqrt{x'^2+y'^2}f(x,y,\theta)$, where $f(x,y,\theta)=\cos\theta E_1(x,y,\theta) + \sin\theta E_2(x,y,\theta)$. Thus, a {$\mathcal{T}$-}{}Lagrangian is characterized by a sum of a classical null Lagrangian and a {parametrization} invariant Lagrangian (c.f., \eqref{eq:para inv lagrangian 2d}). The Euler-Lagrange equations for this case is characterized further by Corollary \ref{appendixcor}.
\end{remark}

\section{\texorpdfstring{{$\mathcal{T}$-}}{T} Lagrangian: second-order Lagrangians}
\label{sec:second order}
In this section, we explore the second-order Lagrangians of the form $\mathscr{L}(s,\uu,\uu',\uu'')$. Such functionals are motivated by the Helfrich-Canham energy \eqref{eq:axi HC} discussed above {in \S\ref{sec:motivation}}.

Following the arguments given in Proposition \ref{lem:pnlagrangian_lem}, we state
\begin{prop}
${\mathscr{L}}(s,\uu,\uu',\uu'')$ is a {$\mathcal{T}$-}{}Lagrangian if and only if the tangential component of its Euler-Lagrange equation vanishes for all $\uu\in \mathcal{C}^{4}_R$. That is,
\begin{equation}
\uu'\cdot\mathfrak{E}_\mathscr{L}(\uu)=\uu'\cdot\Big(D_{ss}\mathscr{L}_{\uu''}-D_s{{\mathscr{L}}}_{\uu'} + {{\mathscr{L}}}_\uu\Big)\equiv 0.
\label{eq:tgtEL2ndOrder=0}
\end{equation}
\label{lem:2nd:pnlagrangian_lem}
\end{prop}

In light of this result and Proposition \ref{lemma2_2}, we infer the following relationship between {parametrization} invariance and {$\mathcal{T}$-}{}Lagrangian:
\begin{cor}
\label{cor:para T lag}
The Lagrangian $\mathscr{L}$ of a {parametrization} invariant functional \eqref{eq:not:functional and lagrangian} is a {$\mathcal{T}$-}{}Lagrangian.
\end{cor}

To proceed with the representation theorem for {$\mathcal{T}$-}{} Lagrangians in the second-order case, {consider the definitions}
\begin{equation}
\mathbf{p}:=\uu'\text{ and }\mathbf{q}:=\mathbf{p}'.
\label{eq:pqdefs}
\end{equation}

We record the following two lemmas:
\begin{lemma}
If ${\mathscr{L}}(s,\uu,\uu',\uu'')$ is a {$\mathcal{T}$-}{} Lagrangian, then
\begin{subequations}
\begin{equation}
(p_i \mathscr{L}_{q_i})_{,q_j} = 0,
\label{eq:q''Coeff}
\end{equation}
\begin{equation}
p_i \Big(\mathscr{L}_{q_ip_j}-\mathscr{L}_{p_iq_j}+2q_k\mathscr{L}_{q_iq_jp_k}\Big) = 0,
\label{eq:q'Coeff}
\end{equation}
\end{subequations}
for each $j\in{\{1,\dotsc,n\}}$, where Einstein summation convention is assumed
on the repeated indices $i$ and $k$.
\end{lemma}
\begin{proof} Refer Appendix \ref{sec:aux claims}.
\end{proof}

\begin{lemma}
If ${\mathscr{L}}(s,\uu,\uu',\uu'')$ is a {$\mathcal{T}$-}{} Lagrangian, then
\begin{equation}
p_i\mathscr{L}_{q_i}=A,
\label{eq:2ndOrder_A}
\end{equation}
\begin{equation}
p_i\mathscr{L}_{p_i}+2q_i\mathscr{L}_{q_i}-\mathscr{L}=q_iA_{,p_i}+B,
\label{eq:2ndOrder_q' simp}
\end{equation}
where {$A=A(s,\uu,\mathbf{p})$ and $B=B(s,\uu,\mathbf{p})$} are (smooth) functions depending on the indicated variables.
\label{lem:5.2}
\end{lemma}
\begin{proof}

Straightforward integration of \eqref{eq:q''Coeff} results in
\eqref{eq:2ndOrder_A}. 

By rewriting the term $p_i\mathscr{L}_{q_ip_j}$ as
$(p_i\mathscr{L}_{q_i})_{,p_j}-\mathscr{L}_{q_i}\delta_{ij}$ and $p_iq_k\mathscr{L}_{q_iq_jp_k}$ as
$q_k\left((p_i\mathscr{L}_{q_i})_{,p_k}-\mathscr{L}_{q_i}\delta_{ik}\right)_{,q_j}$ in
{\eqref{eq:q'Coeff}}, where $\delta_{ij}$ denotes the Kronecker delta, we can express {\eqref{eq:q'Coeff}} as
\begin{equation}
(p_i\mathscr{L}_{q_i})_{,p_j}-\mathscr{L}_{q_i}\delta_{ij}-(p_i\mathscr{L}_{p_iq_j})+2q_k\left((p_i\mathscr{L}_{q_i})_{,p_k}-\mathscr{L}_{q_i}\delta_{ik}\right)_{,q_j}=0.
\notag
\end{equation}

Plugging \eqref{eq:2ndOrder_A} into the previous equation,
we obtain
\[
\left(A_{,p_j}-\mathscr{L}_{q_j}\right)-p_i \mathscr{L}_{p_iq_j} +
2q_k\left(A_{,p_k}-\mathscr{L}_{q_k}\right)_{,q_j}=0.
\]
Since $A$ is independent of
$\mathbf{q}$, we can simplifiy the previous equation as
\[
A_{,p_j}-\mathscr{L}_{q_j}-p_i\mathscr{L}_{p_iq_j} - 2q_k\mathscr{L}_{q_kq_j}=0.
\] 
Writing $q_k\mathscr{L}_{q_kq_j}$ as $(q_k \mathscr{L}_{q_j})_{,q_k}-\mathscr{L}_{q_k}\delta_{jk}$ and shifting terms on the other side, we obtain
\begin{equation}
-A_{,p_j}+\left(p_i\mathscr{L}_{p_i}+2q_i\mathscr{L}_{q_i}-\mathscr{L}\right)_{,q_j}=0,
\notag
\end{equation}
which on integration gives
\begin{equation}
p_i\mathscr{L}_{p_i}+2q_i\mathscr{L}_{q_i}-\mathscr{L}=q_iA_{,p_i}+B,
\notag
\end{equation}
for some function $B=B(s,\uu,\mathbf{p})$ with indicated dependency. This establishes \eqref{eq:2ndOrder_q' simp}.
\end{proof}

{
Consider the (local) change of coordinates:
\begin{equation}
(\mathbf{p},\mathbf{q})\mapsto
(r,\bm{\theta},\xi,\bm{\eta}),
\label{eq:nd:change of variablesform}
\end{equation}
defined by
\begin{subequations}
\begin{equation}
r = \sqrt{\mathbf{p}\cdot \mathbf{p}}, 
\label{eq:r=sq p.p}
\end{equation}
\begin{equation}
\hat{\mathbf{p}}(\bm{\theta}) = \frac{\mathbf{p}}{r},
\label{eq:hatp=p/r}
\end{equation}
\begin{equation}
\xi = \frac{\mathbf{p}\cdot \mathbf{q}}{r^2},
\label{eq:xi def}
\end{equation}
\begin{equation}
(r\hat{\mathbf{p}}_{,\bm{\theta}})[\bm{\eta}] = \left(\mathbf{q}\otimes \hat{\mathbf{p}}(\bm{\theta}) - \hat{\mathbf{p}}(\bm{\theta})\otimes \mathbf{q}\right)[\hat{\mathbf{p}}(\bm{\theta})].
\label{eq:eta def}
\end{equation}
\label{eqs:pq to r t x n}
\end{subequations}
Note that \eqref{eq:hatp=p/r} can be inverted using \eqref{eq:pn equations} to find $\bm{\theta}$ in terms of $\mathbf{p}$. Similarly, the Jacobian matrix $\hat{\mathbf{p}}_{,\bm{\theta}}$ can be inverted to find $\bm{\eta}$ in terms of $\mathbf{p}$ and $\mathbf{q}$. If follows from straightforward differentiation that $\bm{\eta}(s)=\bm{\theta}'(s)$.

The inverse transformation is given by:
\begin{subequations}
\begin{equation}
\mathbf{p} = r \hat{\mathbf{p}}(\bm{\theta}),
\label{eq:p in terms of r theta}
\end{equation}
\begin{equation}
\mathbf{q} = r\hat{\mathbf{p}}\xi + r\hat{\mathbf{p}}_{,\bm{\theta}}\bm{\eta}.
\label{eq:q in terms of xi and eta}
\end{equation}
\label{eqs:inverse general}
\end{subequations}
Under the change of variables \eqref{eq:nd:change of variablesform}--\eqref{eqs:pq to r t x n}, we define
\begin{equation}
\hat{A}(s,\uu,r,\bm{\theta}):=A(s,\uu,\mathbf{p})
\label{eq:nd:hat A}
\end{equation}
\begin{equation}
\hat{B}(s,\uu,r,\bm{\theta}):=B(s,\uu,\mathbf{p}).
\label{eq:nd:hat B}
\end{equation}
}

\begin{remark}
{As an example, when $n=2$, the transformations \eqref{eqs:pq to r t x n} reduce to:}
\begin{equation}
(p_1,p_2,q_1,q_2)\mapsto
(r,\theta,\xi,\eta),
\label{eq:change of variablesform}
\end{equation}
defined by
\begin{subequations}
\begin{equation}
r = \sqrt{p_1^2+p_2^2},
\label{eq:polar r}
\end{equation}
\begin{equation}
\theta = \arctan(p_2/p_1),
\label{eq:polar theta}
\end{equation}
\begin{equation}
\xi = \frac{p_1q_1+p_2q_2}{p_1^2+p_2^2},
\end{equation}
\begin{equation}
\eta = \frac{p_1q_2-p_2q_1}{p_1^2+p_2^2}.
\end{equation}
\label{eq:change of variables}
\end{subequations}
\end{remark}

{
\begin{remark}
Using \eqref{eqs:pq to r t x n} (for $n=2$) and the above Remark, the Lagrangian of this functional (integrand of \eqref{eq:euler}) can be written as $\mathscr{L}_{e}=(\eta/r)^2 r$, which agrees with form suggested by the representation theorem \ref{thm:rep thm second order}.
Note that \eqref{eq:euler} is {parametrization} invariant; an observation has also been made recently in \cite{bates2016elastica}. 
According to Corollary \ref{cor:para T lag} such functionals are {$\mathcal{T}$-}{}Lagrangians. 
\end{remark}
}

In this part, we state and prove the second main result of this paper, the representation theorem for a second-order {$\mathcal{T}$-}{} Lagrangian, i.e.
\begin{theorem}
\label{thm:rep thm second order}
If ${\mathscr{L}}(s,\uu,\uu',\uu'')$ is a {$\mathcal{T}$-} Lagrangian, then ${\mathscr{L}}$ has the local representation
\begin{equation}
{\mathscr{L}}(s,\uu(s),\uu'(s),\uu''(s))=D_s\Xi(s,\uu(s),r(s),{\bm{\theta}}(s)) + r(s) f\Big(\uu(s),\bm{\theta}(s),\frac{{\bm{\eta}}(s)}{r(s)}\Big),
\notag
\end{equation}
where $\Xi$ and $f$ are arbitrary smooth functions.
\label{thmatangt_full 2dof}
\end{theorem}
{A {$\mathcal{T}$-}{}Lagrangian is characterized by a sum of a classical null Lagrangian and a {parametrization} invariant Lagrangian (c.f., \eqref{eq:rep form 2nd order}).}

\begin{proof}
Under the new variables defined by \eqref{eq:nd:change of variablesform}, the
equations \eqref{eq:2ndOrder_A} and \eqref{eq:2ndOrder_q' simp} are respectively
transformed as follows:
\begin{equation}
\mathscr{L}_{\xi} = \hat{A},
\label{eq:2ndOrder-Lagxi}
\end{equation}
\begin{equation}
r\mathscr{L}_r + {\bm{\eta}\cdot \mathscr{L}_{\bm{\eta}}}-\mathscr{L} = \xi(r\hat{A}_{,r}-\hat{A})+{\bm{\eta}\cdot \hat{A}_{,\bm{\theta}}} + \hat{B},
\label{eq:2ndOrder:rL_r+nL_n-L}
\end{equation}
{with $\hat{A}=\hat{A}(s,\uu,r,{\bm{\theta}})$.
Details of this calculation are provided in the Proposition \ref{prop:b:simplify} in the appendix.}

Integrating \eqref{eq:2ndOrder-Lagxi}, we obtain
\begin{equation}
\mathscr{L} = \xi \hat{A}(s,\uu,r,{\bm{\theta}}) + C(s,\uu,r,{\bm{\theta},\bm{\eta}}).
\label{eq:2ndOrder_L=xi_A+B}
\end{equation}

Plugging \eqref{eq:2ndOrder_L=xi_A+B} into \eqref{eq:2ndOrder:rL_r+nL_n-L}, we obtain
\begin{equation}
r C_{,r} + {\bm{\eta}\cdot C_{,\bm{\eta}}} - C = {\bm{\eta} \cdot \hat{A}_{,\bm{\theta}}} + \hat{B},
\label{eq:2ndOrder_rC_rPDE}
\end{equation}
where we omit the variables for brevity.

Solving the first-order partial differential equation \eqref{eq:2ndOrder_rC_rPDE} using the
method of characteristics (see supplementary material), we obtain
\begin{subequations}
\begin{equation}
C(s,\uu,r,{\bm{\theta},\bm{\eta}}) = {\bm{\eta}\cdot \mathcal{A}_{,\bm{\theta}}}(s,\uu,r,{\bm{\theta}}) + \mathcal{B}(s,\uu,r,{\bm{\theta}}) + r E(s,\uu,{\bm{\theta}},\frac{{\bm{\eta}}}{r}),
\label{eq:C=eta A_theta+...}
\end{equation}
where 
\begin{equation}
\mathcal{A}(s,\uu,r,{\bm{\theta}}):=\int \frac{\hat{A}(s,\uu,r,{\bm{\theta}})}{r}\;dr,
\notag
\end{equation}
\begin{equation}
\mathcal{B}(s,\uu,r,{\bm{\theta}}):=\int \frac{\hat{B}(s,\uu,r,{\bm{\theta}})}{r^2}\;dr,
\notag
\end{equation}
{and $E$ is a smooth function.}
\end{subequations}

Note that the integration is performed by treating other variables (i.e.,
$s,\uu,$ and ${\bm{\theta}}$) as constants.

Plugging \eqref{eq:C=eta A_theta+...} into \eqref{eq:2ndOrder_L=xi_A+B}, we obtain
\begin{equation}
\mathscr{L} = \xi \hat{A} + {\bm{\eta}\cdot \mathcal{A}_{,\bm{\theta}}}+\mathcal{B}+r E.
\label{eq:L=xiAetc}
\end{equation}
Differentiating both sides of {\eqref{eqs:pq to r t x n}(a,b)} 
with respect to $s$
results in the following equations: {$\bm{\theta}'(s)=\bm{\eta}(s)$},
$r'(s)=\xi(s)r(s)$. Plugging them into \eqref{eq:L=xiAetc}, we obtain
\begin{equation}
\mathscr{L} = r'\mathcal{A}_{,r}+{\bm{\theta}'\cdot \mathcal{A}_{,\bm{\theta}}} + \mathcal{B}+ r E,
\label{eq:L 2ndorder simplified}
\end{equation}
where, in the first term, we use
\begin{equation}
\frac{1}{r}A=\frac{\partial}{\partial r}\int \frac{A}{r}\;dr = \frac{\partial}{\partial r}\mathcal{A}.
\notag
\end{equation}
{Note that $A=A(s,\uu,r,{\bm{\theta}})$.}

Using the identity (due to the chain rule of differentiation) 
\[
r'\frac{\partial\mathcal{A}}{\partial r}+{\bm{\theta}'\cdot\frac{\partial \mathcal{A}}{\partial \bm{\theta}}}=D_s\mathcal{A}-\frac{\partial \mathcal{A}}{\partial s}-{\mathbf{u}'\cdot \frac{\partial \mathcal{A}}{\partial \mathbf{u}}},
\]
we find that \eqref{eq:L 2ndorder simplified} can be simplified to the form
\begin{equation}
\mathscr{L} = D_s \mathcal{A} (s,\uu,r,{\bm{\theta}})+ Q(s,\uu,r,{\bm{\theta}})+r E(s,\uu,{\bm{\theta}},\frac{{\bm{\eta}}}{r}),
\label{eq:2ndOrder Lag form}
\end{equation}
{
where, using $\mathbf{u}'=\mathbf{p} = r\hat{\mathbf{p}}(\bm{\theta})$, we have defined 
\begin{equation}
Q(s,\uu,r,\bm{\theta}):=-\frac{\partial \mathcal{A}}{\partial s}-r\hat{\mathbf{p}}(\bm{\theta})\cdot \frac{\partial\mathcal{A}}{\partial \mathbf{u}}+\mathcal{B}.
\notag
\end{equation}
}
Plugging \eqref{eq:2ndOrder Lag form} back into \eqref{eq:tgtEL2ndOrder=0}, we obtain
{
\begin{equation}
\bm{\eta}\cdot Q_{,\bm{\theta}} - r\Big(\bm{\eta}\cdot Q_{,r\bm{\theta}}+r\xi Q_{,rr}-\hat{\mathbf{p}}(\bm{\theta})\cdot Q_{,\uu} + \mathbf{p} Q_{,r\uu} + Q_{,sr}\Big)-rE_{,s} = 0.
\label{eq:tgtial lag condition in Q}
\end{equation}
Details of this calculation are provided in Proposition~\ref{prop:b:simplify2} of the appendix.
}

The previous equation must hold for all $s$, $\uu$, and $\p$, and since $Q$ {and E are}  independent of $\xi$, the coefficient of $\xi$ must vanish, i.e.,
\[
Q_{,rr}=0,
\]
whence we obtain the following representation for $Q$:
\begin{equation}
Q(s,\uu,r,{\bm{\theta}}) = r \psi(s,\uu,{\bm{\theta}})+\phi(s,\uu,{\bm{\theta}}),
\label{eq:2ndOrder Q rep}
\end{equation}
for the arbitrary functions $\psi$ and $\phi$. Note that the functions $\psi$ and $\phi$ are independent of $r$, {$\bm{\eta}$}, and $\xi$. 

Using \eqref{eq:2ndOrder Q rep} in \eqref{eq:tgtial lag condition in Q}, the latter can be re-written as
{
\begin{equation}
\frac{\partial}{\partial s}(r E + Q) = \bm{\eta}\cdot \phi_{,\bm{\theta}}+r\hat{\mathbf{p}}(\bm{\theta}) \phi_{,\mathbf{u}}+\phi_{,s};
\notag
\end{equation}
}
from which upon integrating both sides with respect to $s$ (and treating variables $\uu$ and $r$ as constants), we obtain
\begin{align}
& (rE+Q)-(rE+Q)|_{s=0}\notag\\
& = \int_0^s {\Big(\bm{\eta} \cdot \phi_{,\bm{\theta}}(t,\uu,\bm{\theta})} + {r\hat{\mathbf{p}}(\bm{\theta})\phi_{,\uu}(t,\uu,\bm{\theta})\Big)} \;dt+\phi(s,\uu,\bm{\theta})
-\phi(0,\uu,\bm{\theta}).
\label{eq:temp temp}
\end{align}
Since {$\mathbf{u}$}, $r$, and $\theta$ are treated as constants with respect to the integration variable, they may be pulled out of the integral sign. Using {$\uu' = r\hat{\mathbf{p}}(\bm{\theta})$}, the first term in the right side of \eqref{eq:temp temp} can be concisely written as a total derivative:
{
\begin{equation}
(\bm{\eta}\cdot\frac{\partial}{\partial\bm{\theta}} + r \hat{\mathbf{p}}(\bm{\theta})\frac{\partial}{\partial \uu})\int_0^s \phi(t,\uu,\bm{\theta})\;dt = D_s \int_0^s \phi(t,\uu,{\bm{\theta}})\;dt-\phi(s,\uu,{\bm{\theta}}).
\notag
\end{equation}}
Thus, \eqref{eq:temp temp} can be written as
\begin{equation}
rE + Q = D_s\Phi(s,\uu,{\bm{\theta}})+ rE(0,\uu,{\bm{\theta}},\frac{{\bm{\eta}}}{r}) + Q(0,\uu,r,{\bm{\theta}})-\phi(0,\uu,{\bm{\theta}}),
\notag
\end{equation}
where we have defined $\Phi(s,\uu,{\bm{\theta}}) := \int_0^s \phi(t,\uu,{\bm{\theta}})\;dt$. 
Thus, a {$\mathcal{T}$-}{} Lagrangian has necessarily the form
\begin{equation}
\mathscr{L} = D_s\big(\mathcal{A}(s,\uu,r,{\bm{\theta}})+\Phi(s,\uu,{\bm{\theta}})\big) + rE(0,\uu,{\bm{\theta}},\frac{{\bm{\eta}}}{r}) + Q(0,\uu,r,{\bm{\theta}})-\phi(0,\uu,{\bm{\theta}}),
\notag
\end{equation}
which coincides with the stated representation \eqref{thmatangt_full 2dof} by identifying $\Xi$ with $\mathcal{A}+\Phi$ (as function of $(s,\uu,r,{\bm{\theta}})$) and $f(\uu,{\bm{\theta}},\frac{{\bm{\eta}}}{r})$ with $E(0,\uu,{\bm{\theta}},\frac{{\bm{\eta}}}{r})+\psi(0,\uu,{\bm{\theta}})$.
\end{proof}
{
\begin{remark}
It is clear from the above proof that Theorem \ref{thmatangt_full 2dof} is a generalization of known results, for example see Bates et. al \cite{bates2016elastica}, in particular, Theorem A.31 and Remark A.32. Note that \eqref{eq:2ndOrder-Lagxi} and \eqref{eq:2ndOrder:rL_r+nL_n-L} have the same left-hand side (in the coordinates defined by \eqref{eqs:pq to r t x n}) as that of Bates et al. (Theorem A.31 and Remark A.32), but the right-hand side is different. 
\end{remark}
}
\section{\texorpdfstring{{$\mathcal{N}$-}}{N} Lagrangian: case of normal variations}
\label{sec_normalvar}
Let us explore the counterpart to definition \ref{def:T lag} when the first variation trivially vanishes for all normal variations. 
As a departure from the previous discussion, we only focus on the case of first-order Lagrangians in this section. 
Consider
\begin{definition}[{$\mathcal{N}$-}{}Lagrangian]
${\mathscr{L}}$ is an {$\mathcal{N}$-}{}Lagrangian if for any $\uu\in \mathcal{C}_R^2$,
$\delta\mathcal{E}[\delta\uu]=0$ for all $\delta\uu\in N_\uu$.
\end{definition}
Following the proof of Proposition \ref{lem:pnlagrangian_lem}, it is straightforward to establish the 
\begin{prop}
If ${\mathscr{L}}(s,\uu,\uu')$ is a {$\mathcal{N}$-}{}Lagrangian then the normal component of the Euler-Lagrange operator must vanish for all $\uu\in \mathcal{C}_R^2$. That is,
\begin{equation}
\mathbf{n}\cdot\mathfrak{E}_\mathscr{L}(\uu)=\mathbf{n}\cdot\Big(-D_s\mathscr{L}_{\uu'} + {{\mathscr{L}}}_\uu\Big)\equiv 0,
\label{eq:ch2:tgtEL=0}
\end{equation}
for all $\mathbf{n}(s)\in\mathbb{R}^n$ {such that} $\mathbf{n}(s)\cdot\uu'(s)=0$, $s\in[0,1]$. 
\label{thm:ch2:pnlagrangian_thm}
\end{prop}

We prove the counterpart to Theorem \ref{thmatangt_full ndof} to the case of {$\mathcal{N}$-}{}Lagrangians, i.e.
\begin{theorem}
If $\mathscr{L}(s,\uu,\uu')$ is a {$\mathcal{N}$-}{}Lagrangian, then $\mathscr{L}$ has the local representation
\begin{equation}
\mathscr{L}(s,\uu,\uu') = \uu' \cdot {\mathbf{A}}(s,\uu)+ C(s,\uu),
\label{eq:rep lemma eq ndof2}
\end{equation}
where $C$ and {${\mathbf{A}}:= A_i\ee_i$} are arbitrary functions (of variables indicated).
\end{theorem}

\begin{proof}
Since $\mathscr{L}(s,\uu,\uu')$ is a {$\mathcal{N}$-}{}Lagrangian, we have (from Theorem
\eqref{thm:ch2:pnlagrangian_thm}) for any $\uu\in C^2[0,1]$,
\begin{equation}
\Big(-D_s{\mathscr{L}}_{\uu'} + {\mathscr{L}}_\uu\Big)\cdot\a=0,
\notag
\end{equation}
for all $\a\in\mathbb{R}^n$ such that $\a\cdot\uu'=0.$

As before, let us define $\p:=\uu'$, i.e., $\mathscr{L}(s,\uu,\p)$. Expanding $D_s\mathscr{L}_{\p}$ in
the previous equation, we have
\begin{equation}
(-\mathscr{L}_{\p\p}[\p']-\mathscr{L}_{\p\uu}[\p]-\mathscr{L}_{\p s} + \mathscr{L}_{\uu})\cdot\a= 0.
\notag
\end{equation}
Collecting the $\p'$ terms,
\begin{equation}
-\mathscr{L}_{\p\p}[\p']\cdot\a + (\text{ terms that depend on }s, \uu, \p) = 0
\notag
\end{equation}
Since this must hold for all $\uu\in\mathbb{R}^n$, and $\p\in\mathbb{R}^n$, the coefficients
of $\p'$ must vanish:
\begin{equation}
\mathscr{L}_{\p\p}[\a]=0,
\notag
\end{equation}
which, using the definition provided in \eqref{prthetandof}, is equivalent to
\begin{equation}
\mathscr{L}_{\p\p}[{\hat{\p}}_{,\theta_i}({\boldsymbol{\theta}})]=0, i=1, \dotsc, n-1.
\label{eq:ch2:Lpp conditions ndof2}
\end{equation}
By Lemma \ref{lem:consequence of homogeneity2}, we find that
\begin{equation}
\mathscr{L}_{\p} (s,\uu,\p)= {\mathbf{A}}(s,\uu,\sqrt{\p\cdot\p}).
\notag
\end{equation}
Differentiating this with respect to $\p$,
\begin{equation}
\mathscr{L}_{\p\p}(s,\uu,\p) =D_3{\mathbf{A}}(s,\uu,\sqrt{\p\cdot\p})\otimes(\p/\sqrt{\p\cdot\p})=D_3{\mathbf{A}}(s,\uu,\sqrt{\p\cdot\p})\otimes{\hat{\p}}({\boldsymbol{\theta}}),
\notag
\end{equation}
i.e., $\mathscr{L}_{p_ip_j}=(D_3A_i)p_j/\sqrt{\p\cdot\p}.$
Clearly, $(D_3{\mathbf{A}})\otimes{\hat{\p}}[{\hat{\p}}_{,\theta_i}] =0.$ We also require that
\begin{equation}
(D_3{\mathbf{A}})\otimes{\hat{\p}}={\hat{\p}}\otimes (D_3{\mathbf{A}}),
\notag
\end{equation}
i.e. $D_3A_i\hat{p}_j=D_3A_j\hat{p}_i$. 
So $({\hat{\p}}\otimes {\hat{\p}})(D_3{\mathbf{A}})=({\hat{\p}}\otimes (D_3{\mathbf{A}})){\hat{\p}}=((D_3{\mathbf{A}})\otimes{\hat{\p}}){\hat{\p}}=D_3{\mathbf{A}}.$
Hence, $(D_3{\mathbf{A}})\parallel{\hat{\p}}$, i.e.,
\begin{equation}
(D_3{\mathbf{A}})(s,\uu,\sqrt{\p\cdot\p})=A(s,\uu,\p){\hat{\p}}({\boldsymbol{\theta}}).
\notag
\end{equation}
Since ${\boldsymbol{\theta}}$ is arbitrary, we conclude that $A=0$ so that ${\mathbf{A}}$ is independent of $\sqrt{\p\cdot\p}$.
That is,
\begin{equation}
\mathscr{L}_{\p} = {\mathbf{A}}(s,\uu).
\notag
\end{equation}
Integrating this equations,
\begin{equation}
\mathscr{L} = \p\cdot {\mathbf{A}}(s,\uu) + C(s,\uu),
\notag
\end{equation}
for some function $C$.
\end{proof}

We show that above representation theorem does not add anything more to the admissible class of Lagrangians than just the set of null Lagrangians. Thus, the following result recovers the classical divergence representation of null Lagrangians, i.e.
\begin{cor}
If $\mathscr{L}$ is a {$\mathcal{N}$-}{}Lagrangian then it is also a null Lagrangian.
\label{thm:Nnull=nullLag}
\end{cor}
\begin{proof}
Recall that since $\mathscr{L}$ is a {$\mathcal{N}$-}{}Lagrangian, we have
\[
\mathbf{n}\cdot\mathfrak{E}_\mathscr{L}(\uu)=0,\text{ for all }\mathbf{n}\in\mathbb{R}^n,\;{\text{such that}}\;\mathbf{n}\cdot\uu'= 0.
\]
Alternatively,
\begin{equation}
\mathfrak{E}_\mathscr{L}(\uu) - \frac{\p\cdot\mathfrak{E}_\mathscr{L}(\uu)}{|\p|^2}\p=\mathbf{0}.
\label{eq:proj normal}
\end{equation}
Recall that $\p=\uu'\ne0$ as $\uu([0,1])$ is a regular curve.

Due to \eqref{eq:rep lemma eq ndof2}, we can write
\[
\mathfrak{E}_\mathscr{L}(\uu) = (\mathbf{A}_{,\uu})^T[\p] + C_{,\uu}-D_s\mathbf{A},
\]
where, in indicial notation, the $i$-th component of the Euler-Lagrange operator is given by $A_{i,u_j}p_i + C_{,u_i} - D_sA_i$. Expanding the derivative in the last term, we obtain
\[
\mathfrak{E}_\mathscr{L}(\uu) = (\mathbf{A}_{,\uu})^T[\p] + C_{,\uu}-\mathbf{A}_{,s}-\mathbf{A}_{,\uu}[\p] = (\mathbf{A}_{,\uu}^T-\mathbf{A}_{,\uu})[\p] + C_{,\uu}-\mathbf{A}_{,s}.
\]
Plugging this equation into \eqref{eq:proj normal}, we obtain
\[
\Big((\mathbf{A}_{,\uu}^T-\mathbf{A}_{,\uu})[\p] + C_{,\uu}-\mathbf{A}_{,s}\Big)-\frac{\p}{|\p|^2}\Big(\p\cdot(\mathbf{A}_{,\uu}^T-\mathbf{A}_{,\uu})[\p] + \p\cdot(C_{,\uu}-\mathbf{A}_{,s})\Big)=\mathbf{0}.
\]
Let us define
\begin{equation}
\mathbf{B} := \mathbf{A}_{,\uu}^T-\mathbf{A}_{,\uu},\;\mathbf{D} = C_{,\uu}-\mathbf{A}_{,s}.
\label{eq:temp1}
\end{equation}
Writing $\norm{\p}=\sum_k p_k^2$ and taking a common denominator, we obtain, in component notation, 
\[
\sum_{j,k} (B_{ij}p_j+D_i)p_k^2 - p_i\sum_{j,k} (B_{jk}p_jp_k + D_k p_k)=0,\;\text{ for }i=1,\dotsc n.
\]
Note that since $\mathbf{B}$ is skew-symmetric, $\sum_{j,k}B_{jk}p_jp_k=0$. We thus have
\begin{equation}
\sum_{j,k} B_{ij}p_k^2p_j + (D_ip_k-D_kp_i)p_k =0,\;\text{ for } i=1,\dotsc n.
\label{eq:temp}
\end{equation}
Equation~\eqref{eq:temp} must hold for all $\uu\in\mathcal{C}^2_R$, but since $\mathbf{B}$ and $\mathbf{D}$ are independent of $\p(s)=\uu'(s)$, the left hand side of this equation when viewed as a polynomial in $p_i$ must be a zero polynomial, i.e., the coefficients of $p_ip_k$ and $p_k^2p_j$ must vanish:
\begin{equation}
B_{ij}=0,\; D_i=0, \text{ for all }i,j.
\notag
\end{equation}
That is, using \eqref{eq:temp1},
\begin{equation}
\mathbf{A}_{,\uu} = \mathbf{A}_{,\uu}^T,\;\mathbf{C}_{,\uu} = \mathbf{A}_{,s},
\notag
\end{equation}
which in terms of components are given by
\begin{equation}
A_{i,u_j} = A_{j,u_i},\;C_{u_i}=A_{i,s}.
\notag
\end{equation}
Applying Poincar\'e Lemma (i.e. Lemma \ref{lem:4d:poincare ndof}) (when $e_i\equiv 0$), we deduce the existence of $\Psi(s,\uu)$, such that $\mathbf{A}=\Psi_{,\uu}$. Therefore, we have
\[
\mathscr{L}(s,\uu,\uu') = D_s\Psi(s,\uu).
\]
In particular, we obtain the characterization as the classical null Lagrangian. 
\end{proof}

\section{{Discussion}}
\label{sec_discuss}
In this section, we discuss the implications of Theorems \ref{thmatangt_full ndof} and \ref{thm:rep thm second order} to more examples from continuum mechanics.

{\bf Multi-phase fluid membrane}
The Cahn-Hilliard model \cite{cahn1958free} is a popular approach to model phase transitions and spinodal decomposition. In this model, the phase of the material is described by a scalar \emph{order parameter}, a function, $\psi:\omega\to\mathbb{R}$. To model fluid membranes, we assume that the order parameter is defined on the \emph{deformed configuration} of the surface, $\omega$ (c.f., \eqref{eq:axisymmetric para}). The free energy for the system is given by:
\begin{equation}
\mathcal{E}_{CH}(\psi,\mathbf{f}) = \int_\omega \frac{\epsilon}{2}\norm{\nabla \psi}^2 + W(\psi) \;da,
\notag
\end{equation}
where $W=(\psi^2-1)^2$ is a two well function of the order parameter and $\norm{\nabla\psi}^2=g^{\alpha\beta}\psi_{,\alpha}\psi_{,\beta}$ (given in terms contravariant components of the metric tensor, $g_{\alpha\beta}$, of $\omega$). The parameter, $\epsilon$, controls the length scale of the boundary layer between the two phases.

To compare this with our one-dimensional model, we note that 
\[
\norm{\nabla\psi}^2 =\frac{1}{x'^2+y'^2}\psi'^2,
\]
where we assume that the order parameter is a function of axisymmetric surface parametrization, c.f., \eqref{eq:axisymmetric para}. We define 
\[
p_1:=\psi', \quad p_2:=x', \quad p_3:=y'.
\]
In terms of $(r,\bm{\theta})$, c.f., \eqref{eq:pn equations}, we have $\norm{\nabla\psi}^2 = \cot^2\theta_1$ and $da = x r\sin(\theta_1)\;ds$. Thus, the (axisymmetric) Lagrangian corresponding to this functional is given by
\begin{equation}
\mathscr{L}_{CH} = 2\pi r \Bigg(\frac{\epsilon}{2}
\cot^2 \theta_1 + W(\psi) 
\Bigg) x \sin\theta_1.
\label{eq:CH axi}
\end{equation}
Note that $\theta_2=\arctan(y'/x')$ does not appear in the Lagrangian. The combination of Cahn-Hilliard with the Helfrich-Canham energy noted above is a widely studied model for exploring liquid-liquid phase transitions in lipid membranes. In this case, \eqref{eq:CH axi} must be added to the Lagrangian corresponding to the Helfrich-Canham energy:
\begin{equation}
\mathscr{L}_{HC} = 2\pi r\Bigg[\frac{\kappa}{4} \Big(\frac{\eta_2}{r\sin\theta_1}+\frac{\sin\theta_1\sin\theta_2}{x}\Big)^2 + \kappa_g \Big(\frac{\eta_2\sin\theta_2}{x r\sin\theta_1}\Big)\Bigg]x \sin\theta_1. 
\notag
\end{equation}

Since the energy of a fluid membrane is defined with respect to the deformed configuration, it follows that $\mathscr{L}_{CH}+\mathscr{L}_{HC}$ is parametrization invariant. Thus, as a consequence of Corollary \ref{cor:para T lag} such a Lagrangian is a {$\mathcal{T}$-}{}Lagrangian. It then follows that $\mathscr{L}_{CH}+\mathscr{L}_{HC}$ must have the representation given by Theorem~\ref{thmatangt_full 2dof}. Our calculation above is in consonance with this observation. Note the $\eta_2$'s dependence in the Lagrangian only through the ratio $\eta_2/r$ in agreement with Theorem \ref{thm:rep thm second order}. Note that for this model $n=3$, while such results are only known for $n=2$. 

new{The Cahn-Hilliard model is an example of a first-order phase field model, where the energy depends on the first derivative of the order parameter $\psi$. A second-order phase-field model depends on the second derivative of $\psi$. An example of such a theory is the Landau-Brazovskii theory \cite{brazovskii1987theory}, which models crystallization in the bulk and on curved interfaces \cite{dharmavaram2016landau}. According to our representation theorem, c.f, Theorem \ref{thmatangt_full 2dof}, if such a theory is supposed to be parametrization invariant (hence a {$\mathcal{T}$-}{}Lagrangian theory), then the energy density must depend on $\eta_1/r$. While for a general second-order theory, the dependence of the Lagrangian on $\psi''$ can be quite general, parametrization invariance puts restrictions on the form that it can take, which we have characterized in the theorem.}

{\bf Frank-Oseen model} 
Consider a two-dimensional thin film of liquid crystals. We assume the unit director, $\mathbf{n}$, of the liquid crystal molecules is tangential to its surface, we call $\omega$. The \emph{one-constant} Frank-Oseen model for a fluid film is given by \cite{nestler2018orientational}, 
\begin{equation}
\mathcal{E}_{FO}(\mathbf{f},\mathbf{n}) = \frac{K}{2}\int_\omega (\text{div }\mathbf{n})^2 + (\text{rot }\mathbf{n})^2\;da,
\label{eq:frank oseen}
\end{equation}
where $\mathbf{f}$ is the deformation map of the surface $\omega$, $\text{div}$ and $\text{rot }$ are the surface divergence and curl on the deformed surface $\omega$. While, traditionally, the Frank-Oseen free energy is written for a bulk liquid crystal \cite{sheng1975introduction}, there has been a recent interest in generalizing the model for liquid crystal on thin films \cite{nestler2018orientational, borthagaray2021q}. If bending elasticity of the thin fluid film is relevant, the free energy \eqref{eq:frank oseen} can be combined with the Helfrich energy, c.f.,~\eqref{HCenergy} \cite{nitschke2020liquid}.

Specializing to axisymmetry, recall that the deformation map of an axisymmetric surface can be parametrized by \eqref{eq:axisymmetric para}. We parametrized the director field by
\[
\mathbf{n} = \alpha(s)\mathbf{e}_r(\varphi)+ \beta(s)\mathbf{e}_{\varphi} + \gamma(s) \mathbf{k}.
\]
Thus, we envision $\mathcal{E}_{FO}(x(s),y(s),\alpha(s),\beta(s),\gamma)$. Note that the unit director condition means that $\alpha^2+\beta^2+\gamma^2\equiv 1$, which can be enforced using a Lagrange multiplier. It can be shown that (see supplementary materials for details)
\begin{equation}
\text{div}\;\mathbf{n} = \frac{1}{x\sqrt{x'^2+y'^2}} \frac{d}{ds}\Big(\frac{x}{\sqrt{x'^2+y'^2}}(x'\alpha+y'\gamma)\Big),
\notag
\end{equation}
\begin{equation}
\text{rot}\;\mathbf{n} = \frac{(x\beta'+\beta x')}{x\sqrt{x'^2+y'^2}}.
\notag
\end{equation}
To see the relationship between \eqref{eq:frank oseen} and Theorem~\ref{thmatangt_full 2dof}, we define $p_1=\alpha'$, $p_2=\beta'$, $p_3=\gamma'$, $p_4=x'$, $p_5=y'$, and similarly, $q_1=\alpha''$, etc. Thus, this model corresponds to $n=5$.

Expressing in terms of $\bm{\theta}$ and $\bm{\eta}$, c.f., \eqref{eqs:pq to r t x n}, we find that

\begin{equation}
\text{div}\;\mathbf{n} = \frac{1}{r x\sin\theta_1\sin\theta_2\sin \theta_3}\frac{d}{ds}\Big(x(\alpha \cos\theta_4 + \gamma\sin \theta_4)\Big),
\label{eq:divn}
\end{equation}
\begin{equation}
\text{rot }\mathbf{n} = \frac{x \sin\left( \theta_1 \right) \cos\left( \theta_2 \right) + \beta \sin\left( \theta_1 \right) \sin\left( \theta_2 \right) \sin\left( \theta_3 \right) \cos\left( \theta_4 \right)}{x \sin \theta_3 \sin \theta_2 \sin \theta_1}.
\notag
\end{equation}
After expanding the derivative in \eqref{eq:divn} and using the facts that $\theta_4'=\eta_4$ and $da=r\sin\theta_1\sin\theta_2\sin\theta_3$, it is easy to see that the Lagrangian associated with \eqref{eq:frank oseen} has the form of a {parametrization} invariant Lagrangian as expressed by Theorem \eqref{thmatangt_full 2dof}. In particular, the $\bm{\eta}$ dependence of the Lagrangian is in terms of $\eta_4/r$.

\section{{Concluding remarks}}
\label{sec:conclusions}
The analysis presented in the paper indicates the relevance of a certain generalization of the notion of {parametrization} invariance, which is a weaker form. This paper demonstrates this relation in the one-dimensional setting with first- and second-order Lagrangians.

\begin{figure}
\centering
\includegraphics[width=5in]{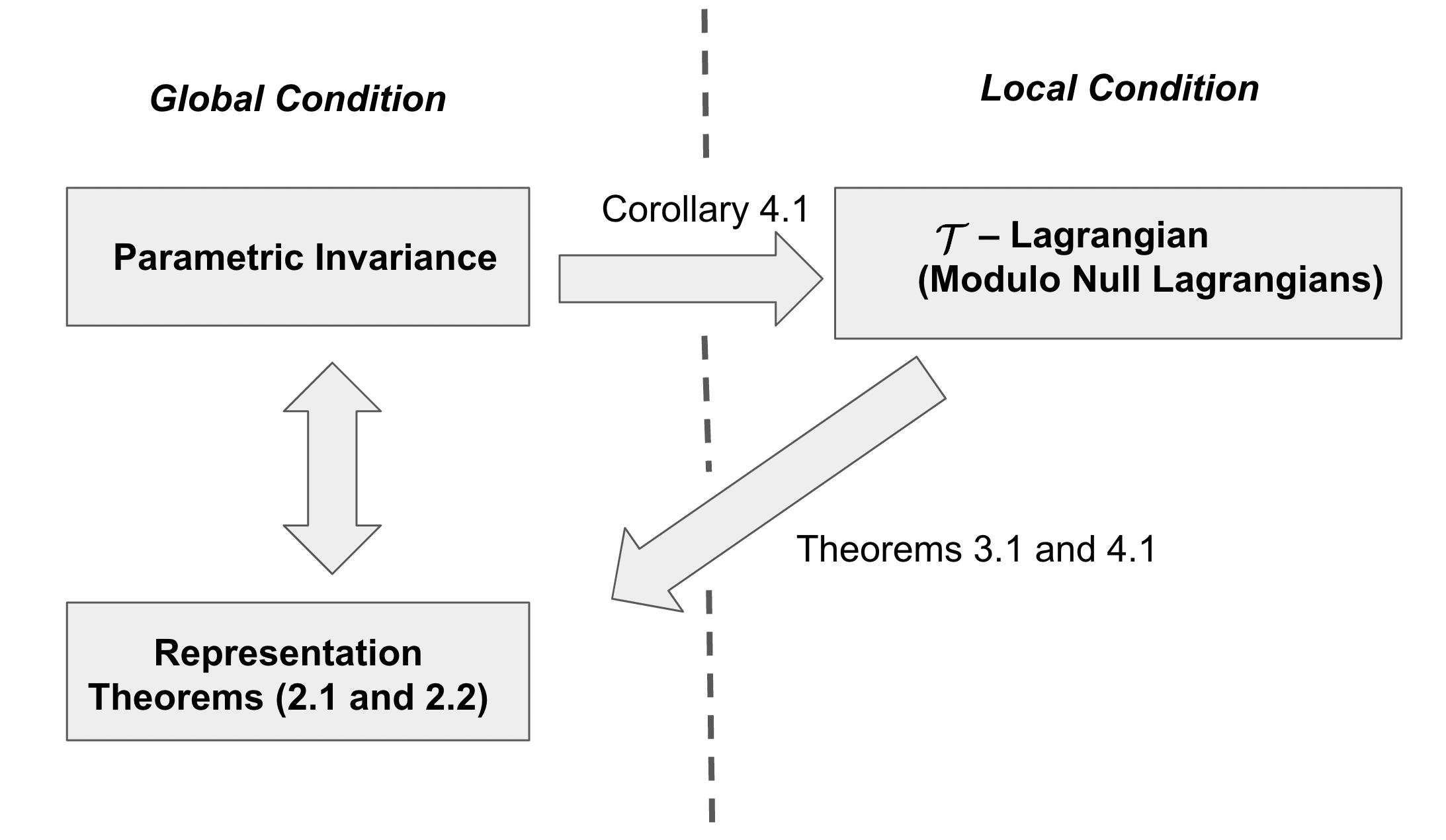}
\caption{Figure showing the equivalence of {parametrization} invariance and {$\mathcal{T}$-}{}Lagrangian (modulo a null Lagrangian).}
\label{fig:significance}
\end{figure}

We introduce a weakened notion of {parametrization} invariance that we term as {$\mathcal{T}$-}{}Lagrangian. While the representation of functionals that possess {parametrization} invariance is well known \cite{GiaquintaV1}, and summarized in Theorems \ref{thm:para inv lagrangian nd} and \ref{thm:2nd:para inv lagrangian nd}, we derive the representation results for {$\mathcal{T}$-}{}Lagrangian. The significance of the weakened notion is that while {parametrization} invariance is a global condition on the functional, {$\mathcal{T}$-}{}Lagrangian is a local condition that can be enforced at the level of equilibrium equations. We see this in Propositions \ref{lem:pnlagrangian_lem} and \ref{lem:2nd:pnlagrangian_lem}, where {$\mathcal{T}$-}{}Lagrangians are equivalently characterized by vanishing of the tangential component of the Euler-Lagrange equations. Through our results, \emph{viz.}, Theorems \ref{thmatangt_full ndof} and \ref{thm:rep thm second order}, we see that on one-dimensional domains, for first and second-order, {parametrization} invariance and {$\mathcal{T}$-}{}Lagrangians (modulo Null Lagrangians) are equivalent. This is summarized schematically in Fig.~\ref{fig:significance}. Thus, a global condition on the Lagrangian (i.e., {parametrization} invariance) can be replaced by a local condition ({$\mathcal{T}$-}{}Lagrangian). The complementary concept of {$\mathcal{N}$-}{}Lagrangian recovers classical null Lagrangians.

new{
We emphasize that while the form of the parametrization invariance of the Lagrangian is known in the case $n=2$, our representation theorem \ref{thmatangt_full 2dof} generalizes this observation for even {$\mathcal{T}$-}{}Lagrangians for $n\geq 1$. 
}

In this paper, we utilize the classical framework of the calculus of variations \cite{Gelfand}; in other words, we deliberately avoid dwelling on the generalization to weakly differentiable functions and measure-theoretic approach \cite{Ball81,Dacorogna}. Besides this, except for the use of Poincar\'e Lemma, we also have avoided the framework of differential geometry and the jet-bundle formalism of variational analysis \cite{saunders1989geometry,Gr99,Crampin}.
These aspects are deemed important in the context of question raised in this paper but are relegated to future investigations.

new{From certain viewpoint, we have proven that in the one-dimensional setting, our way of weakening the notion of parametrization invariance via {$\mathcal{T}$-}{}Lagrangian does not yield any new representations (up to a null Lagrangian).} We conjecture that generalizing the results of the paper to higher dimensions as well as higher order Lagrangians ($p>1, \mathcal{L}(\boldsymbol{u})=\int _{\Omega}{\mathscr{L}}(\boldsymbol{x},\{\nabla^k\boldsymbol{u}(\boldsymbol{x})\}_{k=1}^p)d{x}$) has non-trivial implications and the links with null Lagrangian are also anticipated {\cite{OS88,sharma2021null}}. However, the method of the present paper becomes cumbersome in this regard due to the presence of a high number of indices in the expressions as well as the difficulty of integrating some of the partial differential equations.

\section*{Acknowledgments}
{The authors thank the anonymous reviewers for their constructive comments and suggestions to improve the manuscript.}
BLS gratefully acknowledges the partial support of SERB MATRICS grant MTR/\allowbreak2017/000013.
The part of the work of BLS during Jan-Feb 2023 was partially supported by a grant from the Simons Foundation.
BLS would like to thank the Isaac Newton Institute for Mathematical Sciences, Cambridge, for support and hospitality 
where a part of work on this manuscript was undertaken and was supported by EPSRC grant no EP/R014604/1.
\appendix

\section{Appendix: Auxiliary Claims}
\label{sec:aux claims}
\begin{lemma}
If ${\mathscr{L}}(s,\uu,\uu',\uu'')$ is a {$\mathcal{T}$-}{} Lagrangian, then
\begin{equation}
(p_i \mathscr{L}_{q_i})_{,q_j} = 0,
\label{eq:app:q''Coeff}
\end{equation}
\begin{equation}
p_i \Big(\mathscr{L}_{q_ip_j}-\mathscr{L}_{p_iq_j}+2q_k\mathscr{L}_{q_iq_jp_k}\Big) = 0,
\label{eq:app:q'Coeff}
\end{equation}
for each $j\in\{1,2\}$, where Einstein summation convention is assumed
on the repeated indices $i$ and $k$.
\end{lemma}
\begin{proof}
Writing \eqref{eq:tgtEL2ndOrder=0} explicitly using the Einstein convention, we have
\begin{equation}
p_i\Big( D_s(\mathscr{L}_{q_iq_j}q_j'+\mathscr{L}_{q_ip_j}p_j'+\mathscr{L}_{q_iz_j}z_j'+\mathscr{L}_{q_is}-\mathscr{L}_{p_i}) + \mathscr{L}_{z_i}\Big)=0.
\label{eq:app:tgt cond}
\end{equation}
Expanding the total derivative, $D_s(\cdot)$ and using the notation $z_k'=p_k$, $p_k'=q_k$, we obtain
\begin{multline}
p_i\Bigg(q_j''\mathscr{L}_{q_iq_j} + q_j'\Big(\mathscr{L}_{q_iq_jq_k}q_k' + \mathscr{L}_{q_iq_jp_k}q_k + \mathscr{L}_{q_iq_jz_k}p_k + \mathscr{L}_{q_iq_js}\Big) + \\
\mathscr{L}_{q_ip_j}q_j' + q_j\Big( \mathscr{L}_{q_ip_jq_k}q_k'+\mathscr{L}_{q_ip_jp_k}q_k+\mathscr{L}_{p_ip_jz_k}p_k+\mathscr{L}_{q_ip_js}\Big) + \\
\mathscr{L}_{q_iz_j}q_j + p_j\Big(\mathscr{L}_{q_iz_jq_k}q_k'+\mathscr{L}_{q_iz_jp_k}q_k+\mathscr{L}_{q_jz_jz_k}p_k+\mathscr{L}_{q_iz_js}\Big)+ \\
\Big(\mathscr{L}_{q_isq_k}q_k'+\mathscr{L}_{q_isp_k}q_k+\mathscr{L}_{q_isz_k}p_k+\mathscr{L}_{q_iss}\Big) - \\
\Big(\mathscr{L}_{p_iq_k}q_k'+\mathscr{L}_{p_ip_k}q_k+\mathscr{L}_{p_iz_k}p_k + \mathscr{L}_{p_is}\Big) + \mathscr{L}_{z_i} \Bigg)=0,
\label{eq:app:full expansion}
\end{multline}
where the first three lines in the previous equation arise from the product rule of the first three terms in brackets of \eqref{eq:app:tgt cond}, and the next two lines, by applying $D_s$ to the last three bracketted terms in \eqref{eq:app:tgt cond}.
Since \eqref{eq:app:full expansion} must hold for all $s$, $\uu$, $\p$, and $\q$, the coefficient of $q_{j}''$ and $q_j'$ must be zero. Setting the coefficient of $q_{j}''$ to zero results in
\begin{equation}
p_i \mathscr{L}_{q_iq_j}=0,
\end{equation}
establishing \eqref{eq:app:q''Coeff}. Differentiating the previous equation with respect to $q_k$, $z_k$, and $s$, respectively, results in the following conditions:
\begin{equation}
p_i\mathscr{L}_{q_iq_jq_k}=0,
\end{equation}
\begin{equation}
p_i\mathscr{L}_{q_iq_jz_k}=0,
\end{equation}
\begin{equation}
p_i\mathscr{L}_{q_iq_js}=0,
\end{equation}
where $i,k\in\{1,2\}$. As a consequence of the previous three
equations, certain terms in \eqref{eq:app:full expansion} (which
appear in terms involving $q_j'$) vanish. Setting the
coefficient of $q_j'$ to zero, we obtain \eqref{eq:app:q'Coeff}.
\end{proof}

\begin{cor}
If ${\mathscr{L}}(s,x,y,x',y')$ is a {$\mathcal{T}$-}{}Lagrangian (and therefore has the form \eqref{eq:general form P2d}, where for convenience of notation we write $E_1=-E$ and $E_2=F$), then the critical points $(x(s),y(s))$ of the (normal-component of the) Euler-Lagrange equation are either constant functions or satisfy
\begin{equation}
E_{,y}(x(s),y(s))+F_{,x}(x(s),y(s)) = 0.
\end{equation}
\label{appendixcor}
\end{cor}
\begin{proof} 
Recall \eqref{eq:general form P2d}, i.e., according to Theorem \ref{thmatangt_full ndof}, if ${\mathscr{L}}(s,x,y,x',y')$ is an {$\mathcal{T}$-}{}Lagrangian, then ${\mathscr{L}}$ has the form
\begin{equation}
{\mathscr{L}}(s,x,y,x',y') = \frac{d}{ds}\Psi(s,x(s),y(s)) - \big(x' E(x,y)-y'F(x,y)\big),
\label{eq:rep thm a lag}
\end{equation}
where $\Psi(s,x,y)$, $E(x,y)$, and $F(x,y)$ are arbitrary functions.

Without loss of generality we can ignore the $\Psi$ term as its contribution to the normal component of Euler-Lagrange equation
is trivially zero. So, let 
\begin{equation}
{\mathscr{L}}(s,x,y,x',y') = x' E(x,y)-y'F(x,y).
\end{equation}
The normal component of this Lagrangian is:
\begin{equation}
\begin{split}
&-\big(-(E_{,x} x' + E_{,y} y')+ (x' E_{,x} - y' F_{,x})\big)y' \\
&+ \big((F_{,x} x' + F_{,y} y')+ (x' E_{,y} - y' F_{,y})\big)x' = 0.
\end{split}
\end{equation}
Cancelling $E_{,x}$ and $F_{,y}$ terms, we obtain
\begin{equation}
(E_{,y}+F_{,x})(x'^2+y'^2) = 0,
\end{equation}
i.e.,
\begin{equation}
x'=y'=0\text{ or }(E_{,y}+F_{,x})(x(s),y(s)) = 0.
\end{equation}
\end{proof}

\begin{lemma}
$\p\cdot f_{,\p}(\p)=0$ if and only if $f(\p)=g({\boldsymbol{\theta}})$, for some function $g$, where ${\boldsymbol{\theta}}\in(\mathcal{S}^1)^{n-1}$.
\label{lem:consequence of homogeneity}
\end{lemma}
\begin{proof}

Let us first recall that every vector $\hat{\p}\in\mathcal{S}^{n-1}$ can be associated with a ${\boldsymbol{\theta}}\in{(\mathcal{S}^1)^{n-1}}$, where ${\boldsymbol{\theta}}=(\theta_1, \theta_2, \dotsc, \theta_{n-1})$ in \eqref{eq:pn equations} (schematically shown in Fig.~\ref{Fig2}).

Also note that by Euler's theorem for homogeneous functions, a function $f(\p)$ satisfies the partial differential equation
\begin{equation}
\p\cdot f_{,\p}=0,
\label{eq:ch1:pde1 ndof}
\end{equation}
iff $f(\p)$ is homogeneous in $\p$ of degree $0$.
Consider the change of variables:
\begin{equation}
\p(r,{\boldsymbol{\theta}}) = r{\hat{\p}}({\boldsymbol{\theta}}), \quad {\hat{\p}}({\boldsymbol{\theta}})\cdot{\hat{\p}}({\boldsymbol{\theta}})=1.
\label{prthetandof}
\end{equation}
where, ${\hat{\p}}({\boldsymbol{\theta}})$ is locally described in terms of $n-1$ coordinates $\theta_i\in \mathcal{S}^1$, $i=1, \dotsc, n-1$ noted above. It follows that
\begin{equation}
f_{,r} = {\p}_{,r}[f_{,\p}]= {\p}_{,r}\cdot f_{,\p},\;f_{,{\boldsymbol{\theta}}} = ({\p}_{,{\boldsymbol{\theta}}})^T[f_{,\p}].
\label{frftheta}
\end{equation}
The second equation, in components, is
$f_{,\theta_i}=f_{,p_j}p_{j,\theta_i}.$
Since $\p_{,r} = {\hat{\p}}({\boldsymbol{\theta}})$, $({\p}_{,{\boldsymbol{\theta}}})^T[{\hat{\p}}({\boldsymbol{\theta}})]=\mathbf{0}$ (second equation, in components, holds as $p_kp_k=r^2$ so that $p_{k,\theta_i}p_k=0$), we have
\begin{equation*}
rf_{,r} = f_{,\p}\cdot \p,\; 
f_{,{\boldsymbol{\theta}}}= ({\p}_{,{\boldsymbol{\theta}}})^T[f_{,\p}].
\end{equation*}
Due to \eqref{eq:ch1:pde1 ndof}, we have
\begin{equation*}
rf_{,r}=0 \implies f = g_1({\boldsymbol{\theta}}),\;\text{ for some function }g_1\text{ of }{\boldsymbol{\theta}}.
\end{equation*}
\end{proof}

\begin{lemma}
$\hat{\p}_{,\theta_i}({\boldsymbol{\theta}})\cdot f_{,\p}=0$ for $i=1,\dotsc,n-1$ if and only if $f=g(\sqrt{\p\cdot\p})$, for some function $g$.
\label{lem:consequence of homogeneity2}
\end{lemma}
\begin{proof}
Consider the change of variables:
\begin{equation}
\p(r,{\boldsymbol{\theta}}) = r{\hat{\p}}({\boldsymbol{\theta}}), \quad {\hat{\p}}({\boldsymbol{\theta}})\cdot{\hat{\p}}({\boldsymbol{\theta}})=1.
\label{prthetandof2}
\end{equation}
It follows by differentiating $f$ that
\begin{equation}
f_{,{\boldsymbol{\theta}}} = (\p_{,{\boldsymbol{\theta}}})^T[f_{,\p}].
\end{equation}
That is, for each $i=1,\dotsc n-1$, we have
\begin{equation}
f_{,\theta_i} = r\hat{\p}_{,\theta_i}\cdot f_{,\p}=0,
\end{equation}
where the last equality follows from the hypothesis. Thus, $f$ is independent of ${\boldsymbol{\theta}}$ and must only depend on $r=\sqrt{\p\cdot\p}$.
\end{proof}

{
\section{Claims from Second Order Theory}}
Recall the change of variables given by \eqref{eqs:pq to r t x n} and their inverse \eqref{eqs:inverse general}. Taking the derivatives of $\p$ and $\q$ (as given by the inverse maps \eqref{eqs:inverse general}) with respect to $r$, $\bm{\theta}$, $\xi$, and $\bm{\eta}$, we obtain:
\begin{subequations}
\begin{eqnarray}
\mathbf{p}_{,r} &=& \hat{\mathbf{p}}(\bm{\theta}),
\label{eq:b:p_r}\\
\mathbf{q}_{,r} &=& \xi \hat{\mathbf{p}} + \hat{\mathbf{p}}_{,\bm{\theta}}\bm{\eta},
\label{eq:b:q_r}\\
\mathbf{p}_{,\bm{\theta}} &=& r\hat{\mathbf{p}}_{,\bm{\theta}},
\label{eq:b:p_theta}\\
\mathbf{p}_{,\xi} &=& \mathbf{0},
\label{eq:b:p_xi}\\
\mathbf{q}_{,\xi} &=& r\hat{\mathbf{p}} \quad(= \mathbf{p}),
\label{eq:b:q_xi}\\
\mathbf{p}_{,\bm{\eta}} &=& \mathbf{O},
\label{eq:b:p_eta}\\
\mathbf{q}_{,\bm{\eta}} &=& r \hat{\mathbf{p}}_{,\bm{\theta}},
\label{eq:b:q_eta}
\end{eqnarray}
\end{subequations}
where $\mathbf{O}$ is a zero matrix (of appropriate dimension)

\begin{lemma}
If $\mathscr{L}(s,\uu,\p,\q)$ is a function of the indicated variables, then, in terms of the new variables defined by \eqref{eqs:pq to r t x n}, we have the following identity:
\begin{equation}
\mathbf{p}\cdot \mathscr{L}_{\mathbf{p}} + 2\mathbf{q}\cdot \mathscr{L}_{\mathbf{q}} = r \mathscr{L}_{,r} + \xi\mathscr{L}_{,\xi} + \bm{\eta}\cdot\mathscr{L}_{,\bm{\eta}}
\end{equation}
\label{lem:B:pLp+2qLq}
\end{lemma}
\begin{proof}
{By inspection of \eqref{eq:b:p_r}--\eqref{eq:b:q_eta},} let us first note that $\p$ and $\q$ {(recall \eqref{eqs:pq to r t x n} and \eqref{eqs:inverse general}) both} depend on $r$ and $\bm{\theta}$, while only $\q$ depends on $\xi$ and $\bm{\eta}$.

By the chain rule, it follows that, 
\[
\mathscr{L}_{,r} = \mathbf{p}_{,r}\cdot \mathscr{L}_{,\mathbf{p}} + \mathbf{q}_{,r}\cdot \mathscr{L}_{,\mathbf{q}} = \hat{\mathbf{p}}\cdot \mathscr{L}_{,\mathbf{p}} + \xi \hat{\mathbf{p}}\cdot \mathscr{L}_{,\mathbf{q}} + (\hat{\mathbf{p}}_{,\bm{\theta}}\bm{\eta})\cdot \mathscr{L}_{,\mathbf{q}},
\]
where the second equality follows from \eqref{eq:b:p_r} and \eqref{eq:b:q_r}. Multiplying both sides by $r$ and using $\mathbf{p}=r\hat{\p}$ and \eqref{eq:q in terms of xi and eta} in the last term, we obtain
\begin{equation}
r\mathscr{L}_{,r} = \mathbf{p}\cdot \mathscr{L}_{,\mathbf{p}} + r \xi \hat{\mathbf{p}}\cdot \mathscr{L}_{,\mathbf{q}} + (\mathbf{q}-r\xi\hat{\mathbf{p}})\cdot \mathscr{L}_{,\mathbf{q}}.
\label{eq:b:rLr}
\end{equation}
Similarly, using \eqref{eq:b:p_xi} and \eqref{eq:b:q_xi}, we obtain
\begin{equation}
\xi \mathscr{L}_{,\xi} = \xi \mathbf{p}\cdot \mathscr{L}_{,\mathbf{q}}.
\label{eq:b:xiLxi}
\end{equation}
Finally, differentiating $\mathscr{L}$ with respect to $\bm{\eta}$ and using \eqref{eq:b:p_eta} and \eqref{eq:b:q_eta}, we obtain:
\[
\mathscr{L}_{,\bm{\eta}} = \mathbf{p}_{,\bm{\eta}}^T[\mathscr{L}_{,\mathbf{p}}] + \mathbf{q}_{,\bm{\eta}}^T[\mathscr{L}_{,\mathbf{q}}] = r\hat{\mathbf{p}}^T_{,\bm{\theta}}[\mathscr{L}_{,\mathbf{q}}]
\]
Dotting both sides with $\bm{\eta}$, we obtain 
\begin{equation}
\bm{\eta}\cdot \mathscr{L}_{,\bm{\eta}} = \bm{\eta}\cdot r\hat{\mathbf{p}}^T_{,\bm{\theta}}[\mathscr{L}_{,\mathbf{q}}] = (r\hat{\mathbf{p}}_{,\bm{\theta}}\bm{\eta})\cdot \mathscr{L}_{,\mathbf{q}}=(\mathbf{q}-r\xi\hat{\mathbf{p}})\cdot \mathscr{L}_{,\mathbf{q}},
\label{eq:b:etaLeta}
\end{equation}
where the second equality follows from the definition of the transpose of matrix and the last equality from \eqref{eq:q in terms of xi and eta}. Combining \eqref{eq:b:rLr}, \eqref{eq:b:xiLxi}, and \eqref{eq:b:etaLeta}, we obtain:
\[
r\mathscr{L}_{,r}+\xi \mathscr{L}_{,\xi}+\bm{\eta}\cdot \mathscr{L}_{,\bm{\eta}} = \mathbf{p}\cdot \mathscr{L}_{,\mathbf{p}} + 2 \mathbf{q}\cdot \mathscr{L}_{,\mathbf{q}}.
\]
\end{proof}

The following proposition is a consequence of the previous lemma.
\begin{prop}
In terms of the new variables defined by \eqref{eqs:pq to r t x n}, the Lagrangian must satisfy the condition,
\begin{equation}
r\mathscr{L}_{,r}+\bm{\eta}\cdot\mathscr{L}_{,\bm{\eta}} - \mathscr{L} = \xi(r \hat{A}_{,r}-\hat{A})+\bm{\eta}\cdot \hat{A}_{,\bm{\theta}}+\hat{B},
\end{equation}
where $\hat{A}$ and $\hat{B}$ are defined in \eqref{eq:nd:hat A} and \eqref{eq:nd:hat B}, respectively.
\label{prop:b:simplify}
\end{prop}
\begin{proof}
Note that $A$ and $B$ on the right side of \eqref{eq:2ndOrder_q' simp} are independent of $\q$. That is, their polar forms are independent of $\xi$ and $\bm{\eta}$, \emph{c.f.,} \eqref{eq:nd:hat A}, \eqref{eq:nd:hat B}. It follows by chain rule that,
\[
r \hat{A}_{,r} = \mathbf{p}\cdot A_{,\mathbf{p}}
\]
\[
\hat{A}_{,\bm{\theta}} = r(\hat{\mathbf{p}}^T_{,\bm{\theta}})A_{,\mathbf{p}}
\]
Combining the previous two equations as follows, we obtain
\begin{equation}
\xi r A_{,r} + \bm{\eta}\cdot A_{,\bm{\theta}} = \xi \mathbf{p}\cdot A_{,\mathbf{p}} + (\mathbf{q}-\xi \mathbf{p})\cdot \mathbf{A}_{,\mathbf{p}} = \mathbf{q}\cdot \mathbf{A}_{,\mathbf{p}},
\label{eq:b:rhs of eq}
\end{equation}
where we have used \eqref{eq:q in terms of xi and eta} to replace $\bm{\eta}$ in terms of $\q$. 

Using the Lemma~\ref{lem:B:pLp+2qLq}, the left hand \eqref{eq:2ndOrder_q' simp} of Lemma~\ref{lem:5.2} can be written as $r\mathscr{L}_{,r}+\xi \mathscr{L}_{,\xi}+\bm{\eta}\cdot \mathscr{L}_{,\bm{\eta}}-\mathscr{L}$. Using \eqref{eq:b:rhs of eq} in the right side of \eqref{eq:2ndOrder_q' simp} along with \eqref{eq:2ndOrder-Lagxi}, i.e., $\mathscr{L}_{,\xi}=\hat{A}$, we deduce the require result. 
\end{proof}

The following proposition is used in Theorem \ref{thm:rep thm second order}.

\begin{prop}
If the Lagrangian has the form \eqref{eq:2ndOrder Lag form}, that is,
\[
\mathscr{L} = D_s\mathcal{A} + Q(s,\mathbf{u},r,\bm{\theta}) + r E(s,\mathbf{u},\bm{\theta},\frac{\bm{\eta}}{r}),
\]
then if $\mathscr{L}$ satisfies the {$\mathcal{T}$-}{}Lagrangian condition \eqref{eq:tgtEL2ndOrder=0}, then $Q$ and $E$ must satisfy the condition:
\begin{equation}
\bm{\eta}\cdot Q_{,\bm{\theta}} - r\Big(\bm{\eta}\cdot Q_{,r\bm{\theta}}+r\xi Q_{,rr}-\hat{\mathbf{p}}(\bm{\theta})\cdot Q_{,\uu} + \mathbf{p} Q_{,r\uu} + Q_{,sr}\Big)-rE_{,s} = 0.
\label{eq:b:Qr equation}
\end{equation}
\label{prop:b:simplify2}
\end{prop}
\begin{proof}
By repeated application of product rule for derivatives and using $\q = \p'$, the three terms in \eqref{eq:tgtEL2ndOrder=0} can be individually simplified as follows:
\[
\mathbf{p}\cdot D_{ss}\mathscr{L}_{\mathbf{q}} = D_s(\mathbf{p}\cdot D_s \mathscr{L}_{\mathbf{q}})-\mathbf{q}\cdot D_s\mathscr{L}_\mathbf{q} = D_s(\mathbf{p}\cdot D_s \mathscr{L}_{\mathbf{q}}) - D_s(\mathbf{q}\cdot \mathscr{L}_{\mathbf{q}}) + \mathbf{q}'\cdot \mathscr{L}_{\mathbf{q}}
\]
\[
-\mathbf{p}\cdot D_s \mathscr{L}_\mathbf{p} = -D_s(\mathbf{p}\cdot \mathscr{L}_{\mathbf{p}}) + \mathbf{q}\cdot \mathscr{L}_{\mathbf{p}}
\]
\[
\mathbf{p}\cdot \mathscr{L}_{\mathbf{u}} = D_s\mathscr{L}-\mathscr{L}_s-\mathbf{q}\cdot\mathscr{L}_{\mathbf{p}}-\mathbf{q}'\cdot \mathscr{L}_\mathbf{q},
\]
where the last equation follows from chain rule. Adding the previous three equations, we see that the {$\mathcal{T}$-}-Lagrangian condition \eqref{eq:tgtEL2ndOrder=0} can be written as
\[
D_s\Big(\mathbf{p}\cdot D_s \mathscr{L}_{\mathbf{q}} - \mathbf{q}\cdot\mathscr{L}_\mathbf{q}-\mathbf{p}\cdot \mathscr{L}_{\mathbf{p}}+\mathscr{L}\Big)-\mathscr{L}_s\equiv 0.
\]
Using $\mathbf{p}\cdot D_s \mathscr{L}_\mathbf{q} = D_s(\mathbf{p}\cdot \mathscr{L}_\mathbf{q}) - \mathbf{q}\cdot \mathscr{L}_{\mathbf{q}}$, we obtain
\begin{equation}
D_s\Big(D_s(\mathbf{p}\cdot \mathscr{L}_\mathbf{q})-2 \mathbf{q}\cdot \mathscr{L}_{\mathbf{q}}-\mathbf{p}\cdot \mathscr{L}_\mathbf{p}+\mathscr{L}\Big)-\mathscr{L}_s\equiv 0.
\label{eq:b:null lag temp}
\end{equation}

Since the first term in \eqref{eq:2ndOrder Lag form} is a null Lagrangian, we can neglect it when we plug in the Lagrangian into \eqref{eq:b:null lag temp}.

Note that, neglecting the null Lagrangian term, $\mathscr{L}$ is independent of $\xi$ (c.f., \eqref{eq:2ndOrder Lag form}), i.e., $\mathscr{L}_{,\xi}=0$. Thus, it follows from Lemma~\ref{lem:B:pLp+2qLq} that
\begin{equation}
\mathbf{p}\cdot \mathscr{L}_{\mathbf{p}} + 2 \mathbf{q}\cdot \mathscr{L}_\mathbf{q}-\mathscr{L} = r\mathscr{L}_{,r} + \bm{\eta}\cdot \mathscr{L}_{\bm{\eta}}-\mathscr{L}.
\label{eq:b:pLp+qLq-L}
\end{equation}

Curiously, if $\mathscr{L} = Q + r E$, then,
\begin{equation}
\mathbf{p}\cdot \mathscr{L}_{\mathbf{q}}\equiv 0.
\label{eq:b:pLq}
\end{equation}
To see this first note that
\[
r\hat{\mathbf{p}}(\bm{\theta})\bm{\eta}_{,\mathbf{q}} = \mathbf{I}-\hat{\mathbf{p}}\otimes \hat{\mathbf{p}},
\]
which follows from differentiating \eqref{eq:eta def} with respect to $\q$. Thus, $\bm{\eta}_{,\q}[\p]=\mathbf{0}$. Since $Q$ is independent of $\q$, $\mathscr{L}_\q = \bm{\eta}_{\q}^T[E_{,\bm{\eta}}]$. Thus,
\[
\mathbf{p}\cdot \mathscr{L}_\mathbf{q} = \mathbf{p}\cdot \bm{\eta}_{,\mathbf{q}}^T E_{,\bm{\eta}} = \bm{\eta}_{,\mathbf{q}}[\mathbf{p}]\cdot E_{,\bm{\eta}}=0.
\]

Returning to \eqref{eq:b:null lag temp} along with \eqref{eq:b:pLp+qLq-L}, we can simplify the former equation to
\[
D_s\Big(-r \mathscr{L}_{,r} - \mathbf{\eta}\cdot \mathbf{\mathscr{L}}_{,\mathbf{\eta}}+\mathscr{L}\Big)-\mathscr{L}_s\equiv 0.
\]
Plugging in $\mathscr{L}=Q + rE$ in the the previous equation established the required result \eqref{eq:b:Qr equation}.
\end{proof}

\printbibliography

\end{document}